\documentclass{siamart1116}

\usepackage{amsfonts}
\usepackage{xfrac}
\usepackage{graphicx}
\usepackage{epstopdf}
\usepackage{algorithmic}
\usepackage{amssymb}
\usepackage{booktabs}
\usepackage{algorithmic}
\Crefname{ALC@unique}{Line}{Lines}
\usepackage[T1]{fontenc}
\usepackage{tikz,caption}
\usepackage[caption=false,font=footnotesize]{subfig}
\usepackage{pgfplots}
\pgfplotsset{compat=newest}
\pgfplotsset{plot coordinates/math parser=false}
\usepackage{url}
\usepackage{color}
\usetikzlibrary{plotmarks} 
\usetikzlibrary{external}
\usepgfplotslibrary{external} 
\tikzset{external/system call={lualatex \tikzexternalcheckshellescape -shell-escape -halt-on-error -enable-write18 -interaction=batchmode -jobname "\image" "\texsource"}}

\newtheorem{Lemma}[theorem]{Lemma}

\newlength\figureheight
\newlength\figurewidth 

\ifpdf
  \DeclareGraphicsExtensions{.eps,.pdf,.png,.jpg}
\else
  \DeclareGraphicsExtensions{.eps}
\fi

\def\norm#1{\|#1\|} 
\renewcommand{\vec}[1]{\ensuremath{\mathop{\mathrm{vec}}\left( #1 \right)}}

\newcommand{\diagg}[1]{\ensuremath{\mathop{\mathrm{diag}}\left( #1
\right)}}

\newcommand{\R}{\ensuremath{\mathbb{R}}}

\newcommand{\Am}{\mathcal{A}}

\newcommand{\Ll}{{L}}

\newcommand{\M}{{M}}

\newcommand{\N}{{N}}
\newcommand{\Y}{{Y}}
\newcommand{\U}{{U}}
\newcommand{\V}{{V}}
\newcommand{\Pp}{{P}}

\newcommand{\I}{{I}}
\newcommand{\C}{{C}}
\newcommand{\D}{{D}}

\newcommand{\yct}{\ensuremath{{y}}}
\newcommand{\uct}{\ensuremath{{u}}}
\newcommand{\ubf}{{{u}}}

\newcommand{\yoct}{\ensuremath{{{y}_{\operatorname{obs}}}}}

\newcommand{\W}{{\ensuremath{{W}}}}

\newcommand{\noi}{\ensuremath{\operatorname{noise}}}
\newcommand{\pos}{\ensuremath{\operatorname{post}}}
\newcommand{\pri}{\ensuremath{\operatorname{prior}}}
\newcommand{\Gpri}{\ensuremath{\Gamma_{\operatorname{prior}}}}
\newcommand{\Gnoi}{\ensuremath{\Gamma_{\operatorname{noise}}}}
\newcommand{\Gpos}{\ensuremath{\Gamma_{\operatorname{post}}}}
\newcommand{\pp}{{p}}

\newcommand{\uu}{{u}}
\newcommand{\yy}{{y}}

\newcommand{\nx}{{n_x}}
\newcommand{\nt}{{n_t}}

\def\norm#1{\left\|#1\right\|} 
\numberwithin{theorem}{section}

\newcommand{\TheTitle}{Low-rank computation of posterior covariance matrices in Bayesian inverse problems} 
\newcommand{\ShortTitle}{Low-rank computation in Bayesian inverse problems} 
\newcommand{\TheAuthors}{Peter Benner, Yue Qiu, and~Martin~Stoll}

\headers{\ShortTitle}{\TheAuthors}
\title{{\TheTitle}\thanks{}}

\author{Peter Benner\thanks{Computational Methods in Systems and Control Theory Group,
Max Planck Institute for Dynamics of Complex Technical Systems, Sandtorstr. 1,
39106 Magdeburg, Germany ({\tt benner@mpi-magdeburg.mpg.de})} \and Yue Qiu\thanks{Computational Methods in Systems and Control Theory Group,
Max Planck Institute for Dynamics of Complex Technical Systems, Sandtorstr. 1,
39106 Magdeburg, Germany ({\tt qiu@mpi-magdeburg.mpg.de})}\and Martin Stoll\thanks{Numerical Linear Algebra for Dynamical Systems Group,
Max Planck Institute for Dynamics of Complex Technical Systems, Sandtorstr. 1,
39106 Magdeburg, Germany ({\tt stollm@mpi-magdeburg.mpg.de})}
}

\usepackage{amsopn}

\ifpdf
\hypersetup{
  pdftitle={\TheTitle},
  pdfauthor={\TheAuthors}
}
\fi

\begin{document}

\maketitle

\begin{abstract}
We consider the problem of estimating the uncertainty in statistical inverse problems using Bayesian inference. When the probability density of the noise and the prior are Gaussian, the solution of such a statistical inverse problem is also Gaussian. Therefore, the underlying solution is characterized by the mean and covariance matrix of the posterior probability density. However, the covariance matrix of the posterior probability density is full and large. Hence, the computation of such
a matrix is impossible for large dimensional parameter spaces. It is shown that for many ill-posed problems, the Hessian matrix of the data misfit part has low numerical rank and it is therefore possible to perform a low-rank approach to approximate the posterior covariance matrix. For such a low-rank approximation, one needs to solve a forward partial differential equation (PDE) and the adjoint PDE in both space and time. This in turn gives $\mathcal{O}(\nx\nt)$ complexity for both, computation and storage, where $\nx$ 
is the dimension of the spatial domain and $\nt$ is the dimension of the time domain. Such computations and storage demand are infeasible for large problems. To overcome this obstacle, we develop a new approach that utilizes a recently developed low-rank in time algorithm together with the low-rank Hessian method. We reduce both the computational complexity and storage requirement from $\mathcal{O}(\nx\nt)$ to $\mathcal{O}(\nx + \nt)$. We use numerical experiments to illustrate the advantages of our approach.
\end{abstract}


\pagestyle{myheadings}
\thispagestyle{plain}
\markboth{}{}

\section{Introduction}
Computational mathematicians dealing with simulations of large-scale discretizations describing physical phenomena have made tremendous success over the last decades. This has enabled scientists from various areas of engineering, chemistry, geophysics, et al. to ask more relevant and complex questions. One area that has seen a dramatic increase in the number of published results is the field of statistical inverse problems \cite{Stu10,KaiS05,CalS07}. In particular, the consideration of
partial differential equations (PDEs) as models \textcolor{black}{in statistical inverse problems dramatically increases the problem complexity as a refinement of the model in space and time results in an exponential increase in the problem degrees of freedom}. By this we mean that a discretized problem is typically represented by a spatial system matrix $A\in\R^{n,n}$ where the number of degrees of freedom $\nx$ is typically $O(\frac{1}{h^d})$ with $d$ being the spatial dimension, and
$h$ is the mesh size. It is easily seen that halving the parameter 
$h$ 
means the matrix size will grow by a factor of $2,4,8,\ldots$ depending on the spatial dimension. This complexity is further 
increased when the temporal dimension is incorporated.

While numerical analysis has provided many techniques that allow the efficient handling of such problems, e.g. Krylov methods \cite{book::saad}, multigrid techniques \cite{hackbusch:mm}, we are faced with an even steeper challenge when uncertainty in the parameters of the model is incorporated. For this we consider the approach of Bayesian inverse problems where the goal is to use prior knowledge to obtain information about the conditional mean or the posterior covariance given a set of measured data.
While computing the conditional mean typically corresponds to a problem formulation frequently encountered in PDE-constrained optimization \cite{book::FT2010,HerK10,borzi2003grid}, the problem of computing the posterior covariance matrix is much more challenging as this matrix is dense and involves the inverse of high-dimensional discretized PDE problems. In \cite{FlaWAHBWG11} and subsequent works, the authors proposed a low-rank approximation of the posterior covariance matrix. While this already reduces the complexity dramatically, the storage requirements for the resulting vectors still suffer from high dimensionality with respect to $n_x$ and $n_t$. \textcolor{black}{Our aim in this paper is to introduce a method based on low-rank techniques \cite{StoB15} that allows to reduce the complexity from $O(\nx \nt)$ to $O(\nx+\nt)$.}

In order to do this, we will first derive the basic problem following~\cite{FlaWAHBWG11}. This will be followed by the presentation of a low-rank technique that we previously introduced for PDE-constrained optimization. We then use this to establish a low-rank eigenvalue method based on the classical Lanczos procedure \cite{lanczos} or Arnoldi procedure \cite{Saa11}. After introducing different choices of covariance matrices, we show that our approach can be theoretically justified. We
then illustrate the applicability of our proposed methodology to a diffusion problem and a convection diffusion problem, and present numerical results illustrating the performance of the scheme.

\section{Statistical/Bayesian Inverse Problems}\label{sec_sta_inv}
We refer to \cite{KaiS05,Stu10} for excellent introductions into the subject of statistical inverse problems. We follow \cite{FlaWAHBWG11,BuiGMS13} in the derivation of our model setup and start with the \textit{parameter-to-observable map} $g:\R^{n}\times\R^{k}\rightarrow\mathbb{R}^{m}$ defined as
\begin{equation}
Y=g(U,E),
\end{equation}
where $U,Y,E$ are vectors of random variables. Note that here, $u\in\R^{n}$, our vector of model parameters to be recovered, is a realization of $U,$ the error vector $e\in\R^{k}$ is a realization of $E,$ and $y\in\R^{m}$ is a realization of the vector of observables $Y$. The vector $\yoct\in\R^{m}$ contains the observed values.  As discussed in \cite{BuiGMS13}, even when using the `true' model parameters $u,$ the observables $\yct$ will differ from the measurements $\yoct$ due to measurement noise and the inadequacy of the underlying PDE model.

In a typical application such as the one discussed later, evaluating $g(.)$ requires the solution of a PDE potentially coupled to an observation operator representing a domain of interest. 

The Bayes' theorem, which plays a key role in the Bayesian inference, is written as
\begin{equation}
 \pi_{\pos}:=\pi(u|\yoct)=\frac{\pi_{\pri}(u)\pi(\yoct|u)}{\pi(\yoct)}\propto\pi_{\pri}(u)\pi(\yoct|u),
\end{equation}
where we used the prior probability density function (PDF) $\pi_{\pri}(x)$, the likelihood function $\pi(\yoct|u)$, and the data $\yoct$ with $\pi(\yoct)>0$. The function $\pi_{\pos}:\R^{n}\rightarrow \R$ is the posterior probability density function. The likelihood is derived under the assumption of additive noise
\begin{equation}
 \label{modelf1}
Y=f(U)+E
\end{equation}
where $f:\R^{n}\rightarrow\R^{m}$ and $E$ is the additive noise and given as $E=Y-f(U).$ We once more follow \cite{FlaWAHBWG11}, assuming that $U$ and $E$ are statistically independent and we can use
\[
\pi_{\noi}(e)=\pi_{\noi}(\yoct-f(u)).
\]
Therefore, Bayes' theorem can be written as
\begin{equation}
 \label{bayes1}
\pi_{\pos}\propto \pi_{\pri}(u)\pi_{\noi}(\yoct-f(u)).
\end{equation}
Assuming that both probability density functions for $U$ and $E$ are Gaussian, we can rewrite the PDFs in the form
\begin{equation}
\begin{split}
\pi_{\pri}(u)\propto& \exp\left(-\frac{1}{2}\left(u-\bar{u}_{\pri}\right)^{T}\Gpri^{-1}\left(u-\bar{u}_{\pri}\right)\right),\\
\pi_{\noi}(e)\propto& \exp\left(-\frac{1}{2}\left(e-\bar{e}\right)^{T}\Gnoi^{-1}\left(e-\bar{e}\right)\right),
\end{split}
\end{equation}
where $\bar{u}_{\pri}\in\R^{n}$ is the mean of the model parameter prior PDF and $\bar{e}$ is the mean of the noise PDF. We further have the two covariance matrices 
$\Gpri\in\R^{n,n}$ for the prior and $\Gnoi\in\R^{m,m}$ for the noise. The Gaussian assumption allows us to rewrite Bayes' theorem further to get
\begin{equation}
\label{bayes2}
\begin{split}
\pi_{\pos}\propto& \exp\left(-\frac{1}{2}\left(u-\bar{u}_{\pri}\right)^{T}\Gpri^{-1}\left(u-\bar{u}_{\pri}\right)-\frac{1}{2}\left(e-\bar{e}\right)^{T}\Gnoi^{-1}\left(e-\bar{e}\right)\right)\\
=&\exp\left(-\frac{1}{2}\norm{u-\bar{u}_{\pri}}_{\Gpri^{-1}}^{2}-\frac{1}{2}\norm{e-\bar{e}}_{\Gnoi^{-1}}^{2}\right).
\end{split}
\end{equation}
Let us further assume that the parameter-to-observable map $g(U,Y)$ is given as in \eqref{modelf1} with $f(U)=AU$. The matrix $A\in\R^{m,n}$ represents a linear map from the parameters $u$ to the observables $y$. We will later see that often this matrix \textcolor{black}{involves the inverse of a discretized representation of a PDE operator. In \cite{spantini2016goal}, the authors incorporate a quantity of interest operator into $A$}. Therefore, it will typically be dense and very large. We arrive now at a restated version of Bayes theorem \eqref{modelf1}
\begin{equation}
\label{bayes3}
\begin{split}
\pi_{\pos}\propto&\exp\left(-\frac{1}{2}\norm{u-\bar{u}_{\pri}}_{\Gpri^{-1}}^{2}-\frac{1}{2}\norm{\yoct-Au-\bar{e}}_{\Gnoi^{-1}}^{2}\right).
\end{split}
\end{equation}
From this relation we can express several relevant statistical quantities. For example, we can compute the \textit{maximum a posteriori point} (MAP), which is defined via 
\begin{equation}
 \label{map1}
 \bar{u}_{\pos}=\mathrm{arg max}_{u}\pi_{\pos}(u)
\end{equation}
and to compute it, one can solve the following optimization problem
\[
 \bar{u}_{\pos}=\mathrm{argmin}_{u}\left(\frac{1}{2}\norm{u-\bar{u}_{\pri}}_{\Gpri^{-1}}^{2}+\frac{1}{2}\norm{\yoct-Au-\bar{e}}_{\Gnoi^{-1}}^{2}\right).
\]
Note that this problem is a deterministic inverse problem and resembles the structure one finds in PDE-constrained optimization problems \cite{book::FT2010,book::IK08}. Many strategies are known how to solve this efficiently and in particular underlying function space considerations help in the development of efficient numerical algorithms. An infinite-dimensional discussion of the above problem is given in \cite{BuiGMS13,Stu10} and we only refer to the infinite-dimensional setup when needed.
Our goal in this paper will not be the solution of the MAP problem. The goal of devising low-rank methods for this case has recently been established in \cite{StoB15}  and the techniques there are likely to be applicable as the only difference is the use of the weighting matrices 
$\Gnoi\textrm{ and }{\Gpri},
$
which are for the classical PDE-constrained optimization problem mass matrices or matrices involving mass matrices. 
The more challenging question lies in the approximation of the posterior covariance matrix
\begin{equation}
 \label{postcov1}
 \Gpos=\left(A^{T}\Gnoi^{-1}A+\Gpri^{-1}\right)^{-1}.
\end{equation}


The approximation of $\Gpos$ is in general very costly and, without further approximation, intractable.  The approach presented in \cite{FlaWAHBWG11,BuiGMS13} computes a low-rank approximation to this matrix using the following relation
\begin{equation}
\begin{split}
\Gpos&=\left(A^{T}\Gnoi^{-1}A+\Gpri^{-1}\right)^{-1}\\
&=\Gpri^{1/2}\left(\Gpri^{1/2}A^{T}\Gnoi^{-1}A\Gpri^{1/2}+I\right)^{-1}\Gpri^{1/2}.
\end{split}
\end{equation}
The authors in~\cite{FlaWAHBWG11,BuiGMS13} then compute a low-rank approximation to the so-called \textit{prior-preconditioned} Hessian of the data misfit 
$\tilde{\mathcal{H}}_{\text{mis}}\in\R^{n,n}$
\begin{equation}\label{eqn::prior_pre}
\tilde{\mathcal{H}}_{\text{mis}}=\Gpri^{1/2}A^{T}\Gnoi^{-1}A\Gpri^{1/2}
\end{equation}
with the approximation
\[
\tilde{\mathcal{H}}_{\text{mis}}\approx V\Lambda V^{T},
\]
where $V$ and $\Lambda$ represent the dominant eigenvectors and eigenvalues, respectively. Using this approximation and the Sherman-Morrison-Woodbury formula \cite{book::golubvanloan} one obtains for the prior-preconditioned system 
\[
\left(\Gpri^{1/2}A^{T}\Gnoi^{-1}A\Gpri^{1/2}+I\right)^{-1}\approx\left(V\Lambda V^{T}+I\right)^{-1}=I-V\tilde{\Lambda}V^{T}
\]
with $\tilde{\Lambda}=\mathrm{diag}(\frac{\lambda_i}{\lambda_i+1})$, \textcolor{black}{and $\lambda_i$ is the $i$-th diagonal entry of $\Lambda$}.
A similar approach for 
\begin{equation}
\begin{split}
\left(A^{T}\Gnoi^{-1}A+\Gpri^{-1}\right)^{-1}&\approx\left(V\Lambda V^{T}+\Gpri\right)^{-1}\\
&=\Gpri-\Gpri V\left(\Lambda^{-1}+V\Gpri V^{T}\right)^{-1}\tilde{\Lambda}V^{T}\Gpri,
\end{split}
\end{equation}
becomes a feasible alternative if $\left(\Lambda^{-1}+V\Gpri V^{T}\right)^{-1}$ can be evaluated in reasonable time.

Our main goal of this paper will be the derivation of an efficient scheme to approximate the matrix $V$ from the low-rank approximation to the misfit Hessian. Before discussing this problem, we introduce an idea that becomes instrumental in realizing this and is motivated by a PDE-constrained optimization problem.

\section{A Low-Rank Technique for PDE-Constrained Optimization}\label{sec::modprob}
In order to better understand the stochastic inverse problem, we investigate it in relation to a PDE-constrained optimization problem. We start the derivation of the low-rank in time method by considering an often
used model problem in PDE-constrained optimization (see
\cite{book::hpuu09,book::IK08,book::FT2010}), minimization of
\begin{equation}
\label{MF1}
\min_{\yct,\uct}~\frac{1}{2}\left\|\yct-\yoct\right\|_\mathcal{Q}^{2}
+\frac{\beta}{2}\left\|\uct\right\|_\mathcal{P}^{2}, 
\end{equation}
with $\mathcal{P}$ and $\mathcal{Q}$ space-time cylinders. The constraint linking state $\yct$ and control $\uct$ of this problem is given by the heat
equation with a distributed control term
\begin{equation}
\begin{split}
\label{HEQdc} \yct_t-\nabla^{2}\yct={}&\uct,\quad\text{in }\Omega, \\
\ \nonumber \yct={}&f,\quad\text{on }\partial\Omega.
\end{split}
\end{equation}
Here $\Omega$ denotes the domain and $\partial\Omega$ corresponds to the boundary of the domain. For a more detailed discussion on the well-posedness, existence of solutions, discretization 
et al., we refer the interested reader to \cite{book::hpuu09,book::IK08,book::FT2010}. 
The solution of such an optimization problem is obtained using a Lagrangian approach and considering the first order conditions, which for our problem results in a linear system of the form 
\begin{equation}
 \label{KKT2}
 \begin{split}
\underbrace{
\left[
\begin{array}{ccc}
\D_1\otimes\tau\M_1&0&-\left(\I_{n_t}\otimes\Ll+\C^{T}\otimes \M\right)\\
0&\D_2\otimes\beta\tau\M_2&\D_3\otimes\tau\N^{T}\\
-\left(\I_{n_t}\otimes\Ll+\C\otimes \M\right)&\D_3\otimes\tau\N&0\\
\end{array}
\right]
}_{\mathcal{A}}
 \textcolor{black}{
 \left[
\begin{array}{c}
\mathbf{\yy}\\
\mathbf{\uu}\\
\mathbf{\pp} \\
\end{array}
\right]} &\\
 =
\left[
\begin{array}{c}
\D_1\otimes\tau\M_1\mathbf{y_{\text{obs}}}\\
\mathbf{0}\\
\mathbf{d}\\
\end{array}
\right]&,
\end{split}
\end{equation}
where $\D_1=\D_2=\diagg{1,1,\ldots,1,1}$ and
$\D_3=\I_{n_t}$ come from the discretization of the temporal parts of the objective function or the right hand side of the PDE-constraint (cf.~\cite{BHT10,PSW11,S11_TDTP}). The matrices $M_1$ and $M_2$ are mass matrices corresponding to observation and control domain. The matrix $N$ is essentially representing the incorporation of the control into the constraint, i.e., $N$ is a mass matrix in the above example. The matrix $C$ represents the all-at-once discretization of the time-derivative in the PDE and $L$ the discretized Laplacian. 
 Here, the state, control, and adjoint state are represented by the following space-time vectors
$$
\mathbf{\yy}=
\left[\begin{array}{c}
\mathbf{\yy}_1\\
\vdots\\
\mathbf{\yy}_{n_t}\\
\end{array}
\right],
\mathbf{\uu}=
\left[\begin{array}{c}
\mathbf{\uu}_1\\
\vdots\\
\mathbf{\uu}_{n_t}\\
\end{array}
\right],\textnormal{ and }
\mathbf{\pp}=
\left[\begin{array}{c}
\mathbf{\pp}_1\\
\vdots\\
\mathbf{\pp}_{n_t}\\
\end{array}
\right].
$$

We point out again that the Kronecker product is defined as
\[
W\otimes V=
\left[
\begin{array}{ccc}
w_{11}V&\hdots&w_{1m}V\\
\vdots&\ddots&\vdots\\
w_{n1}V&\hdots&w_{nm}V\\
\end{array}
\right]
\]
and remind the reader of the definition of the $\vec{}$ operator via
\[
\vec{W}=
\left[
\begin{array}{c}
w_{11}\\
\vdots\\
w_{n1}\\
\vdots\\
w_{nm}\\
\end{array}
\right]
\]
as well as the relation 
\[
\left(W^T\otimes V\right)\vec{Y}=\vec{VYW}.
\]
In~\cite{StoB15}, it was shown that the solution to the PDE-constrained optimization problem can be computed in low-rank form
\begin{equation}
\begin{split}
Y=&
\begin{bmatrix}
\mathbf{\yy}_1, \mathbf{\yy}_2, \ldots, \mathbf{\yy}_{n_t}\\
\end{bmatrix}
\approx\W_{\Y}\V_{\Y}^{T}\textnormal{ with }
\W_{\Y}\in\R^{n_1,k_1},\V_{\Y}\in\R^{n_t,k_1},\\
U=&
\begin{bmatrix}
\mathbf{\uu}_1, \mathbf{\uu}_2, \ldots, \mathbf{\uu}_{n_t}\\
\end{bmatrix}
\approx \W_{\U}\V_{\U}^{T}\textnormal{ with }
\W_{\U}\in\R^{n_2,k_2},\V_{\U}\in\R^{n_t,k_2},\\
P=&
\begin{bmatrix}
\mathbf{\pp}_1, \mathbf{\pp}_2, \ldots, \mathbf{\pp}_{n_t}\\
\end{bmatrix}
\approx\W_{\Pp}\V_{\Pp}^{T}\textnormal{ with }
\W_{\Pp}\in\R^{n_1,k_3},\V_{\Pp}\in\R^{n_t,k_3},
\end{split}
\end{equation}
where the $k_i$ are small in comparison to the spatial and temporal dimensions. The authors in \cite{StoB15} illustrated that the low-rank structure of a right hand side is maintained throughout a Krylov subspace iteration and the above described representation. Low-rank techniques for Krylov-subspace methods have recently received much attention and we refer the reader to \cite{KreT10,KreT10a,AndT12} and for tensor structured equations \cite{larskres-survey-2013,DoOs-dmrg-solve-2011,dc-ttgmres-2013,DolS14}.

We obtain a significant storage reduction if we can base our approximation of the
solution using such low-rank factors. It is easily seen that due to the low-rank nature of the factors, we have to perform fewer multiplications with
the submatrices by also maintaining smaller storage requirements.

There are several similarities of the problem \eqref{KKT2} and the statistical inverse problem presented earlier. It is clear that with the choice 
\begin{equation}
\label{CovPDE}
\Gpri=\left(\D_2\otimes\beta\tau\M_2\right)^{-1}\textrm{ and }\Gnoi=\left(\D_1\otimes\tau\M_1\right)^{-1},
\end{equation}
the PDE-constrained optimization problem can be interpreted as a statistical inverse problem and the posterior covariance matrix $\Gpos$ is given by eliminating both state and adjoint state from the system matrix \eqref{KKT2} to obtain a reduced Hessian system. Furthermore, it is clear that for many choices of prior and noise covariance, we can utilize the tensor structure to compute low-rank solutions. For this we state the posterior covariance matrix of the PDE optimization problem

\begin{equation}
\begin{split}
\Gpos=&
\left[
\left(\D_2\otimes\beta\tau\M_2\right)+
\left(\D_3\otimes\tau\N^{T}\right)
\left(\I_{n_t}\otimes\Ll+\C^{T}\otimes \M\right)^{-1}\right.\\
&\left.\left(\D_1\otimes\tau\M_1\right)
\left(\I_{n_t}\otimes\Ll+\C\otimes \M\right)^{-1}\left(\D_3\otimes\tau\N\right)
\right]^{-1}
\end{split}
\end{equation}
and the misfit Hessian
\begin{equation}
\begin{split}
\tilde{\mathcal{H}}_{\text{mis}}=&
\left[
\left(\D_2\otimes\beta\tau\M_2\right)^{1/2}\left(\D_3\otimes\tau\N^{T}\right)
\left(\I_{n_t}\otimes\Ll+\C^{T}\otimes \M\right)^{-1}\right.\\
&\left.\left(\D_1\otimes\tau\M_1\right)
\left(\I_{n_t}\otimes\Ll+\C\otimes \M\right)^{-1}\left(\D_3\otimes\tau\N\right)\left(\D_2\otimes\beta\tau\M_2\right)^{1/2}
\right].
\end{split}
\end{equation}

We keep this example in mind when we now discuss a low-rank technique to approximate the eigenvectors of the posterior covariance matrix in low-rank form. For this we propose a low-rank \textcolor{black}{Krylov subspace method to compute the dominating eigenvectors and eigenvalues.}

Before discussing the eigenvalue approximation strategy, we want to comment on the scaling of the PDE-constrained optimization problem in relation to the statistical inverse problem discussed in \cite{FlaWAHBWG11}. The authors there consider 
\begin{equation}
\begin{split}
\label{MF2}
\min_{\yct,\uct}~\frac{\beta_{\noi}}{2}\left\|\yct-\yoct\right\|_\mathcal{Q}^{2}
+\frac{\beta_{\pri}}{2}\left\|\uct\right\|_\mathcal{P}^{2}, 
\end{split}
\end{equation}
which in the simple case of full observations and control leads to the following rescaling of \eqref{CovPDE}
\begin{equation}
\Gpri=\left(\D_2\otimes\beta_{\pri}\tau\M_2\right)^{-1}\textrm{ and }\Gnoi=\left(\D_1\otimes\tau\beta_{\noi}\M_1\right)^{-1}.
\end{equation}
Assuming that $M\approx h^{d}I$ and $D_i=I$ we get
\begin{equation}
\Gpri^{-1}=\beta_{\pri}\tau h^{d}I\textrm{ and }\Gnoi^{-1}=\beta_{\noi}\tau h^dI.
\end{equation}
In PDE-constrained optimization one typically reduces $\beta$ in~\eqref{MF1} to allow for a more expensive control that drives the state closer to the desired state. This would mean that in the statistical inverse setting $\beta_{\noi}=1$ and decreasing the value of $\beta_{\pri},$ which implies that in \eqref{postcov1} the role of the prior covariance gets diminished and most contributions are coming from the noise. We have similar settings and observations for stochastic inverse problems in this manuscript. These will be shown by the analysis in \Cref{sec::eigenfun} and numerical experiments in \Cref{sec::experiments}. The right choice of parameters $\beta_{\text{prior}}$ and $\beta_{\text{noise}}$ depends on the underlying application and we refer to~\cite{CalS07} for a discussion of the roles of the parameters $\beta_{\pri}$ and $\beta_{\noi}$ as regularization parameters.

\section{Low-rank Lanczos/Arnoldi Method}\label{sec_lanczos}
We recall that our goal is to find a low-rank approximation of the eigenvectors of the posterior covariance matrix. The goal is to compute an approximation to $\tilde{\mathcal{H}}_{\text{mis}}\approx V\Lambda V^{T}$ with $V=[v_1,v_2,\ldots,v_{k}]$ and $k$ much smaller than the dimension of $\tilde{\mathcal{H}}_{\text{mis}}.$ For the PDE-constrained optimization problem $\tilde{\mathcal{H}}_{\text{mis}}\in\R^{\nx\nt,\nx\nt}$. 

Our main assumption at this point is that storing each $v_j$ and especially a number of such space-time vectors can pose serious problems. Additionally, in order to perform the matrix vector multiplication with $\tilde{\mathcal{H}}_{\text{mis}}$, a large number of PDE-solutions need to be computed. For this we point out  that in order to apply the matrix $\tilde{\mathcal{H}}_{\text{mis}}$ in an Arnoldi procedure, we need to solve the spatial system over the whole time-domain. A major advantage of our approach is motived by the fact that
\[
v_j=\vec{V_j}\quad\forall j=1,\ldots,k \textrm{ with } V_{j}\in\R^{\nx,\nt},
\]
which we assume is well approximated via
\begin{equation}\label{eqn::lr_vec}
V_j\approx W_{j,1}W_{j,2}^{T}
\end{equation}
with $W_{j,1}\in\R^{\nx,r_j}$ $W_{j,2}\in\R^{\nt,r_j}$ with $r_j\ll\min\left\lbrace \nx,\nt\right\rbrace.$  If the Arnoldi- or Lanczos vectors are of this form, then the application of the matrix $\tilde{\mathcal{H}}_{\text{mis}}$ to such vectors requires fewer PDE solves than in the full case.

\textcolor{black}{Note that we need to compute the dominant eigenvectors of the prior-precondi\-tioned data misfit Hessian $\tilde{\mathcal{H}}_{\text{mis}}$~\eqref{eqn::prior_pre}, therefore the Lanczos method can be used. At the $j$-th Lanczos iteration, we need to perform $\tilde{\mathcal{H}}_{\text{mis}}v_{j-1}$ using the low-rank approach to get a form like ~\eqref{eqn::lr_vec}. Here $v_{j-1}$ is the $(j-1)$-th Lanczos vector. Due to the low-rank approximation, the orthogonality of Lanczos vectors is lost. Reorthogonalization should be used to orthogonalize Lanczos vectors. Meanwhile, the symmetric property of $\tilde{\mathcal{H}}_{\text{mis}}$ cannot be preserved for the low-rank form of the matrix-vector product. Therefore, we make use of the more general Arnoldi method to compute the dominant eigenvectors of $\tilde{\mathcal{H}}_{\text{mis}}$. We also observe that when applying the Arnoldi method with the truncation error appropriately chosen we still get real eigenvalues. This 
will be shown in \Cref{sec::experiments}.}

We now briefly recall the Arnoldi method, which is the more general procedure. We refer to \cite{book::golubvanloan} for details. We recall that the Arnoldi process for a matrix $B$ can be written as
\[
BV_k=V_{k+1}H_{k+1,k},
\]
where $V_k$ consists of orthonormal columns and $H_{k+1,k}\in\mathbb{R}^{(k+1)\times k}$ is a Hessenberg matrix. The iterative build-up of the columns of $V$ is captured by the recursion
\begin{equation}\label{lanczos1}
\tilde{v}_{k+1}=Bv_{k}-\sum_{i=1}^{k}h_{i,k}v_i,
\end{equation}
where $h_{i,k}=v_i^{T}Av_k$. The vector $\tilde{v}_{k+1}$ is then normalized using the scalar $h_{k+1,k}.$ While this is well-known our goal here is to illustrate how this method is amenable to the use within a low-rank framework. For the Arnoldi process considered in this manuscript, $B=\tilde{\mathcal{H}}_{\text{mis}}$ and the application of $B$ to $v_k$ results in a low-rank matrix, i.e., 
\[
Bv_{k}=\vec{W_{1,B}W_{2,B}^{T}}\]
with small rank. \textcolor{black}{This is because $B$ is related with the inverse of a PDE operator in space and time, which has shown in~\cite{StoB15} that applying such an operator to a low-rank vector again gives a low-rank vector}. We can then write the right hand side of \eqref{lanczos1} as
\begin{equation}
\label{lanczos2}
\vec{W_{1,B}W_{2,B}^{T}}-\alpha_k \vec{W_{1,k}W_{2,k}^{T}}-\beta_{k}\vec{W_{1,k-1}W_{2,k-1}^{T}},
\end{equation}
and write the last expression as
\[
\vec{\left[W_{1,B},-\alpha_kW_{1,k},-\beta_{k} W_{k-1,B}\right]\left[W_{2,B},W_{2,k},W_{2,{k-1}}\right]^{T}}.
\]
The size of the matrix 
\[
\left[W_{1,B},-\alpha_kW_{1,k},-\beta_{k} W_{k-1,B}\right]\in\R^{\nx,r_B+r_k+r_{k-1}}
\]
is increased to $r_B+r_k+r_{k-1}.$  Using truncation techniques this can typically be controlled. For example, one could achieve the truncation by utilizing skinny QR factorizations \cite{KreT10} or truncated singular value decompositions \cite{StoB15}. 

In our numerical experiments, we show that maintaining the orthogonality of the Lanczos vectors is crucial. \textcolor{black}{Meanwhile, since we use low-rank methods to approximate the Lanzcos vectors, the orthogonality and the symmetry of the Lanzcos process get lost}. We hence opt for the more expensive but stable Arnoldi process, where we orthogonalize with respect to all previous vectors. \textcolor{black}{The full reorthogalization also demands more storage than in the Lanczos case}. This is another advantage of our approach. \textcolor{black}{With full reorthogonalization, the storage cost are increasing for both full and low-rank scheme but in the low-rank framework stay significantly below the full scheme}. Here, we give the low-rank Arnoldi method in~\Cref{alg::lr_arnoldi}. 

\begin{algorithm}
	\caption{Low-Rank Arnoldi Method}\label{alg::lr_arnoldi}
	\begin{algorithmic}[1]
		\STATE {\textbf{Input:} maximal Arnoldi steps $m_a$, unit vector $v_1$, truncation tolerance $\varepsilon_0$}
		\FOR {$j=1:\ m_a$}
		\STATE{perform low-rank matrix vector product $w=\tilde{\mathcal{H}}_{\text{mis}}v_j$ up to the truncation tolerance  $\varepsilon_0$}
		\FOR{i=1:\ j} 
		\STATE {perform low-rank dot product $H_{i,\ j}=w^H v_i$} 
		\STATE {update $w\leftarrow w-H_{i, j}v_i$}
		\ENDFOR
		\STATE {$H_{j+1,\ j}=\sqrt{w^H w}$}
		\IF{$j\ <\ m_a$}
		\STATE {$v_{j+1} =\sfrac{1}{H_{j+1,\ j}}w$}
		\ENDIF
		\ENDFOR
		\STATE {\textbf{Output:} low-rank Arnoldi vectors $v_j$, and Hessenberg matrix $H$}
		\end{algorithmic}
\end{algorithm}

\textcolor{black}{We note that the biggest challenge for the low-rank Arnoldi method is to perform the low-rank matrix vector product in line $3$ of~\Cref{alg::lr_arnoldi} since $\tilde{\mathcal{H}}_{\text{mis}}$ is large and dense. We propose the tensor-train (TT) format in~\Cref{sec::experiments} to perform such computations efficiently. The full orthogonalization procedure in line $4$-$7$ of~\Cref{alg::lr_arnoldi} is also performed with the TT format.}

\textcolor{black}{We use the standard Arnoldi method for low-rank eigenvector computations in~\Cref{alg::lr_arnoldi}. This is practical for the problems studied in this manuscript since we just need to compute up to a few hundred Arnoldi vectors. Since the computational complexity of full orthogonalization increases with the number of Arnoldi vectors, if more Arnoldi vectors are needed, the restarted Arnoldi method can be implemented with a low-rank version~\cite{book::golubvanloan}.}



\section{Analysis of the Eigenfunctions}\label{sec::eigenfun}
The eigenfunction analysis for the general case presented above is not straightforward. Our goal in this section is to give a theoretical justification for simple cases. We start with the case of a steady state problem involving the two-dimensional Poisson equation. For this we consider the misfit Hessian
\[
\tilde{\mathcal{H}}_{\text{mis}}=\Gpri^{1/2}A^{T}\Gnoi^{-1}A\Gpri^{1/2}.
\]
Assuming for now that $\Gpri=\beta_{\pri} I_d$ and $\Gnoi=\beta_{\noi} I_d$ so that we are left with
\[
\tilde{\mathcal{H}}_{\text{mis}}=\frac{\beta_{\pri}}{\beta_{\noi}}A^{T}A.
\]
Note that we have assumed for the data misfit Hessian to be defined in function space with the prior and noise covariance operators to be multiples of the identity operator. Our goal is to find eigenfunctions of $A^{T}A$ with $A$ the inverse of the \textcolor{black}{PDE operator} and $A^{T}$ the inverse of the adjoint \textcolor{black}{PDE operator}. For this we use the eigenfunctions of the two-dimensional Laplacian operator satisfying
\[
-\Delta y=\lambda y
\]
with zero Dirichlet boundary conditions (see \cite{KutS84}). \textcolor{black}{Note that for this problem $A^{T}=A$.}  We here consider the eigenvalue problem on a rectangular domain $[0,a]\times[0,b]$ and state the eigenfunctions as 
\[
y_{m,n}=\sin\left(\frac{m\pi x_1}{a}\right)\sin\left(\frac{n\pi x_2}{b}\right)
\]
with $x_1\in[0,a]$ and $x_2\in[0,b]$.
The associated eigenvalues are given by
\[
\lambda_{m,n}=\pi^{2}\left[\left(\frac{m}{a}\right)^2+\left(\frac{n}{b}\right)^2\right]
\]
with $m,n=1,2,\ldots$ (see \cite{KutS84}). Coming back to the misfit Hessian we can write this as
\[
\tilde{\mathcal{H}}_{\text{mis}}=\frac{\beta_{\pri}}{\beta_{\noi}}A^2
\]
with $A$ the inverse Laplacian. Assuming zero Dirichlet boundary conditions it holds that
\[
\tilde{\mathcal{H}}_{\text{mis}}y_{m,n}=\mu_{m,n} y_{m,n}
\]
is satisfied for the eigenfunctions of the Laplacian 
\[
y_{m,n}(x_1,x_2)=\sin\left(\frac{m\pi x_1}{a}\right)\sin\left(\frac{n\pi x_2}{b}\right)
\]
with the eigenvalues given by
\[
\mu_{m,n}=\frac{\beta_{\pri}}{\beta_{\noi}}\lambda_{m,n}^{-2},
\]
where the decay of $\lambda_{m,n}^{-2}$ is quite rapid. This justifies the approximation of $\tilde{\mathcal{H}}_{\text{mis}}$ by a small number of eigenfunctions. In \cite{FlaWAHBWG11}, the authors also justify this for a one-dimensional convection diffusion equation. 

The goal of our approach is to add even more low-rankness to the computation. We want to illustrate that the eigenfunctions used to represent the misfit Hessian are of small rank. When studying the eigenfunctions above it is clear that the eigenfunction
\[
y_{m,n}(x_1,x_2)=\sin\left(\frac{m\pi x_1}{a}\right)\sin\left(\frac{n\pi x_2}{b}\right)
\]
is already separated into $x_1$ and $x_2$ components with \textcolor{black}{separation rank $1$ (cf. \cite{Hack12}), where it is the sum of two functions with one depending on the first and the other on the second variable}. For this case we have established the following lemma
\begin{Lemma}
The eigenfunctions of the misfit Hessian 
\[
\tilde{\mathcal{H}}_{\text{mis}}=\frac{\beta_{\pri}}{\beta_{\noi}}A^2
\]
with $A$ the inverse Laplacian with zero Dirichlet conditions defined on the rectangle $[0,a]\times[0,b]$, are separated and given by 
\[
y_{m,n}(x_1,x_2)=\sin\left(\frac{m\pi x_1}{a}\right)\sin\left(\frac{n\pi x_2}{b}\right)
\]
and hence are of separation rank $1$.
\end{Lemma}
It is not so straightforward to establish similar results for more complicated equations.

For the space-time PDE-constrained optimization problem discussed earlier, we note that
\[
VDV^{T}\approx \alpha(h,\tau,\beta)
\left(\I_{\nt}\otimes\Ll+\C^{T}\otimes \I\right)^{-1}\left(\I_{\nt}\otimes\Ll+\C\otimes \I \right)^{-1},
\]
where we used $\M\approx h^{d}\I$ and collected all scalars in $\alpha(h,\tau,\beta)$. Our aim is to establish eigenvalue and eigenvector results for 
\[
\left(\I_{\nt}\otimes\Ll+\C\otimes \I\right)\left(\I_{\nt}\otimes\Ll+\C^{T}\otimes \I\right).
\]
We note that this fits the well known relation that the singular values of a matrix $\Am\in\R^{m,m}$ are the square roots of the eigenvalues of the matrix $\Am^{T}\Am,$ which, assuming full rank of $\Am$, is a symmetric and positive definite matrix. Now assuming the SVDs 
\[
C=U_{C}\Sigma_CV^{T}_{C}\textnormal{ and }\Ll=U_{\Ll}\Sigma_\Ll V^{T}_{\Ll},
\]
we obtain  
\[
\underbrace{\left(\I_{\nt}\otimes U_{\Ll}\Sigma_\Ll V^{T}_{\Ll}+V_{C}\Sigma_CU^{T}_{C}\otimes\I\right)}_{\Am}=\underbrace{\left(U_{\Ll}\otimes V_{C}\right)}_{U}\underbrace{\left(\I_{\nt}\otimes\Sigma_\Ll +\Sigma_C\otimes \I\right)}_{\Sigma}\underbrace{\left(V^{T}_{\Ll}\otimes U_C^{T}\right)}_{V^{T}}.
\]
From this it follows that 
\[
\Am^{T}\Am=V\Sigma U^{T}U\Sigma V^{T}=V\Sigma^2 V^{T}
\]
is the eigendecomposition of $\Am^{T}\Am,$ which has the same eigenvectors as $\Am^{-1}\Am^{-T}.$ As the eigenvalues of $\Am^{-1}\Am^{-T}$ quickly decay to zero because of the compactness of the operator, we only need a small number of columns of $V$. Our aim in this paper is to express each column of $V$ further in a low-rank fashion. For this we note that $e_1^{(\nx\nt)}=e_1^{(\nx)}\otimes e_1^{(\nt)}$ and ignoring superindices we get
\[
Ve_1=\left(V^{T}_{\Ll}\otimes U_C^{T}\right)\left(e_1\otimes e_1\right)=v^{T}_{1,\Ll}\otimes u_{1,C}^{T},
\]
and hence a vector of rank one if the eigenvectors are all real. Complex eigenvectors would further introduce a small rank increase and the consideration of $\M$ instead of $h^d\I$ can with a simultaneous diagonalization of the pencil $(\Ll,\M)$ lead to small eigenvector ranks of the overall system. 
This justifies our choice of approximating the eigenvectors in low-rank form.

\section{Numerical Results}\label{sec::experiments}
In this section, we study the performance of the low-rank Lanczos algorithm presented in~\Cref{sec_lanczos}. The results presented in this section are based on an implementation of the above described algorithms within MATLAB\textsuperscript{\textregistered}. We perform the discretization of the PDE-operators within the IFISS~\cite{ifiss} framework using $Q1$ finite elements for the heat equation and the streamline upwind Petrov–Galerkin  (SUPG) method for the convection diffusion equation. Our experiments are performed for a final time $T=1$ with a
varying number of time-steps. As the domain $\Omega$ we consider the unit square but other domains are of course possible. We specify the boundary conditions for each problem separately. Throughout the results section we fixed the truncation at $10^{-8}$ for which we observed good results. Additionally, we also performed not listed experiments with a tolerance of $10^{-10}$ for which we also observed fast convergence. Larger tolerances should be combined with a deeper
analysis of the algorithms and a combination with flexible outer solvers. All results are performed on a standard Ubuntu desktop Linux machine with Intel(R) Core(TM) i7-4790 CPU @ 3.60GHz and 8GB of RAM.

The mathematical model we consider in this section is given by the following instationary PDE,

\begin{equation}\label{eqn::pde}
    \begin{split}
        \frac{\partial}{\partial t}y + \mathcal{L} y & = 0, \qquad \qquad \Omega \times (0, T) \\
                                                   y & = u, \qquad \qquad \Omega \times \{t=0 \} \\
                                                   y & = 0, \qquad \qquad \partial \Omega_D \times (0, T) \\
                        \triangledown y \cdot {\bf n} & = 0, \qquad \qquad \partial \Omega_N \times (0, T) 
    \end{split}
\end{equation}
where $\mathcal{L}$ is a PDE operator and for the numerical experiments we consider the case of the heat equation $\mathcal{L}=-\Delta$ and the convection-diffusion equation $\mathcal{L}=-\mu \Delta + \overrightarrow{\omega} \cdot \triangledown$. For all the numerical tests, the initial concentration $u$ represents the unknown parameter $\mathbf{u}$, and the observation data $\mathbf{y}_{\text{obs}}$ are collected by sensors, which are distributed in part of the domain $\Omega$. 

As stated in~\Cref{sec_sta_inv}, the statistical inverse problem with Gaussian noise and prior using a Bayesian formulation is related to a weighted least squares problem. Here we use the same functional as in~\cite{FlaWAHBWG11}, which is given by the following functional

\begin{equation}\label{eqn::obj}
\min_{u} \left(\frac{\beta_{\text{noise}}}{2}
\int_{0}^{T}\int_{\Omega}\left(y-y_{\text{obs}}\right)^2b(x, t)dxdt +\frac{\beta_{\text{prior}}}{2}\int_{\Omega}u^2dx\right),
\end{equation}
in which $u$ satisfies the PDE model~\cref{eqn::pde}, $b(x, t)$ is the observation operator, and $u$ is the uncertainty term, which is the initial condition for this numerical example. Here, we study the sparse observation case, where $b(x, t)$ is defined by 
\[
    b(x, t) = \sum_j \delta_j.
\]
\textcolor{black}{Here $\delta_j$ is the regional function of sensors, $\delta_j=1$ at the region of the $j$-th sensor and $\delta_j=0$ elsewhere.}

Discretizing the functional~\cref{eqn::obj} in turn gives 
\begin{equation}\label{eqn::dis_obj}
\min_{u} \left(\frac{1}{2}
\left({\bf y} - {\bf y}_{\text{obs}}\right)^T\mathcal{B}^T\Gamma_{\text{noise}}^{-1}\mathcal{B}({\bf y} - {\bf y}_{\text{obs}}) +\frac{1}{2} \ubf^T\Gamma_{\text{prior}}^{-1}\ubf\right),
\end{equation}
where $\mathcal{B}$ is the discretization of $b(x, t)$, the variable ${\bf y}$ is the discrete variant of $y$ in~\cref{eqn::pde} stacked in time, and it satisfies $\mathcal{K} {\bf y}=\mathcal{C} u$. Here $\mathcal{K}$ and $\mathcal{C}$ come from the discretization of the PDE model~\cref{eqn::pde}. 

Discretization of the objective function gives that $\Gamma_{\text{noise}}=\sfrac{1}{\beta_{\text{noise}}}I_{n_t}\otimes M$, and $\Gamma_{\text{prior}}=\sfrac{1}{\beta_{\text{prior}}}M$, where $M$ is the mass matrix, $I_{n_t}$ is an $n_t\times n_t$ identity matrix. Here, $n_t$ is the number of time variables. Since we need to compute $\Gamma_{\text{prior}}^{\sfrac{1}{2}}$, we use  $\Gamma_{\text{prior}}\approx\sfrac{1}{\beta_{\text{prior}}}h^{-d}I_{n_x}$. Here, $h$ is the mesh size, $d$ is the spatial dimension, $I_{n_x}$ is a $n_x\times n_x$ identity matrix, and $n_x$ is the number of spatial variables. Other settings such as smoothing operators used for $\Gamma_{\text{prior}}$ have been studied in~\cite{BuiGMS13} and can be used with our approach as well.

As analyzed in~\Cref{sec_sta_inv}, the posterior covariance matrix $\Gamma_{\text{post}}$ is given by the inverse of the Hessian of~\cref{eqn::dis_obj}. Therefore, 
\begin{equation}\label{eqn::postvar}
    \Gamma_{\text{post}} = \left(
    \mathcal{C}^T\mathcal{K}^{-T}\mathcal{B}^T\Gamma_{\text{noise}}^{-1}\mathcal{B}\mathcal{K}^{-1}\mathcal{C} + \Gamma_{\text{prior}}^{-1}
    \right)^{-1}
\end{equation}

Note that for different PDE models, $\mathcal{K}$ is also different after discretization. In this section, we use two types of PDE models, i.e., the heat equation and the convection-diffusion equation to show the performance of our low-rank algorithm for the approximation of $\Gamma_{\text{post}}$. We argue that our low-rank algorithm also applies to other time-dependent PDE operators. Here we apply our method to symmetric systems and unsymmetric systems.

\textcolor{black}{We also point out that due to the uncertainty being given as the initial condition, the posterior Hessian is only of spatial dimension. Nevertheless, the Arnoldi process applied to this matrix requires the solution of space-time problems and the low-rank form of our approach results in a much reduced number of spatial solves. The complexity reduction is even more pronounced when the uncertainty is part of the system as a space-time variable such as the right hand side of the PDE.}
\subsection{Implementation Details}
According to~\eqref{eqn::postvar}, the prior-preconditioned data misfit part after discretization of~\eqref{eqn::obj} is given by 
\begin{equation}\label{eqn::misfit}
\tilde{\mathcal{H}}_{\text{mis}}= \Gamma_{\text{prior}}^{\sfrac{1}{2}}\mathcal{C}^T\mathcal{K}^{-T}\mathcal{B}^T\Gamma_{\text{noise}}^{-1}\mathcal{B}\mathcal{K}^{-1}\mathcal{C}\Gamma_{\text{prior}}^{\sfrac{1}{2}}.
\end{equation}
To apply the Lanczos iteration to~\eqref{eqn::misfit}, we need to solve the space-time discretized PDE $\mathcal{K}$ and adjoint PDE $\mathcal{K}^T$. Here we take $\mathcal{K}$ as an example. This asks us to solve a linear system of the following type,
\begin{equation*}
\left(I_{n_t}\otimes L + C\otimes M\right) x = f,
\end{equation*}
or 
\begin{equation}\label{eqn::st_pde}
\left(I_{n_t}\otimes L + C\otimes M\right) vec(X) = vec(F).
\end{equation}
Here $X$ and $F$ are matrices/tensors of appropriate sizes. Numerical solutions of~\eqref{eqn::st_pde} have been studied intensively from the matrix equation point of view, c.f.~\cite{simoncini16, messweb,BenS13} for an overview.

In this manuscript, we solve~\eqref{eqn::st_pde} based on its tensor structure and use the alternative minimal energy (AMEn) approach~\cite{DolS14} implemented in the tensor-train (TT) toolbox~\cite{tt-toolbox} to solve the tensor equation~\eqref{eqn::st_pde}. At each AMEn iteration, either a left Galerkin projection or a right Galerkin projection is applied to the system~\eqref{eqn::st_pde}. Therefore, after Galerkin projection, we need to solve a linear system either of the format
\begin{equation}\label{eqn::amen_1}
\left(\hat{I}_n\otimes L + \hat{C}\otimes M\right)x = \tilde{b},
\end{equation}
or
\begin{equation}\label{eqn::amen_2}
\left(I_n\otimes \tilde{L} + C\otimes \tilde{M}\right)\tilde{x}=\hat{b}.
\end{equation}
Here $\hat{I}$, $\hat{C}$, $\tilde{L}$, $\tilde{M}$ are matrices of appropriate dimensions after Galerkin projection (c.f.~\cite{DolS14} for details). 

\textcolor{black}{After Galerkin projection, the size of the system~\eqref{eqn::amen_1} is still relatively large while the size of~\eqref{eqn::amen_2} is quite moderate. Therefore, Krylov solvers such as the generalized minimized residual (GMRES)~\cite{book::saad} method or the induced dimension reduction IDR(s)~\cite{sonneveld2008} method can be used to solve~\eqref{eqn::amen_1} while a direct method can be used to solve~\eqref{eqn::amen_2}.} 
	
To accelerate the convergence of the Krylov solver, we use the following preconditioner
\begin{equation}\label{eqn::pre_1}
P=\text{diag}(\hat{I}_n)\otimes L + \text{diag}(C)\otimes M,
\end{equation}
to solve~\eqref{eqn::amen_1}. Here $\text{diag}( \cdot )$ is an operation that extracts the diagonal entries of a matrix and form a diagonal matrix. One can use standard techniques such as multigrid methods~\cite{wesseling1992} or incomplete LU factorization (ilu)~\cite{book::golubvanloan} to approximate the preconditioner~\eqref{eqn::pre_1}. Here we use backslash implemented in MATLAB\textsuperscript{\textregistered}.

\textcolor{black}{We also want to point out that there are many other methods to efficiently solve~\eqref{eqn::st_pde}, such as the low-rank factored alternating directions implicit (ADI) method (cf. \cite{BenK14}).}

\subsection{The Heat Equation}
\begin{figure}[H]
    \centering
    \includegraphics[width=0.5\textwidth]{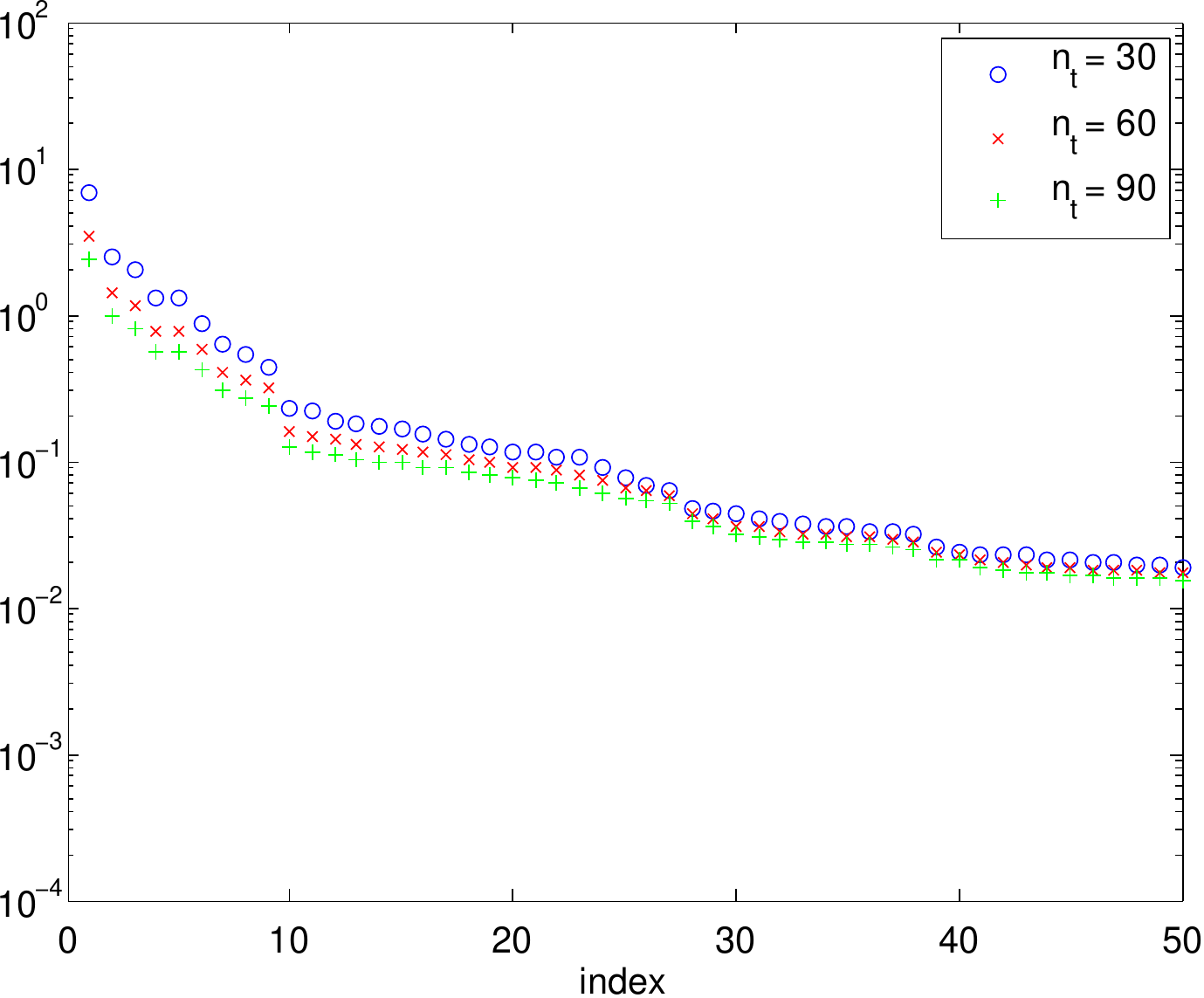}
    \caption{Largest $50$ eigenvalues of $\tilde{\mathcal{H}}_{\text{mis}}$ computed using low-rank Arnoldi}\label{fig::heat_eig_nt369}
\end{figure}

In this part, we use the 2D time-dependent heat equation in a unit square as an example to study the performance of our low-rank algorithm for the heat equation. Discretizing the equation in space using $Q_1$ finite elements and in time using the implicit Euler method gives us an $n_x \times n_t$ linear system, where $n_x$ is the number of spatial variables while $n_t$ is the number of time steps. First we study spectral properties of the prior-preconditioned data misfit part $\tilde{\mathcal{H}}_{\text{mis}}$ and the posterior covariance matrix. Using a
$64\times 64$ grid to discretize the heat equation and set $n_t=30,\ 60,\ 90$, respectively, we plot the $50$ largest eigenvalues of $\tilde{\mathcal{H}}_{\text{mis}}$ in~\cref{fig::heat_eig_nt369}. Here $\beta_{\text{noise}}=10^4\beta_{\text{prior}}$.

As shown in~\cref{fig::heat_eig_nt369}, \textcolor{black}{there are only a few dominant eigenvalues of $\tilde{\mathcal{H}}_{\text{mis}}$. For most cases, a threshold of $10^{-1}$, or even $10^{0}$ is acceptable to approximate $\tilde{\mathcal{H}}_{\text{mis}}$ and to reduce the uncertainty of the system, which will be shown later.} Meanwhile, the number of time steps does not influence the decay rate of $\tilde{\mathcal{H}}_{\text{mis}}$. This makes it possible to compute a fixed number of Arnoldi vectors for even
more time steps. 

\begin{figure}[htb!]
    \centering
    \includegraphics[width=0.5\textwidth]{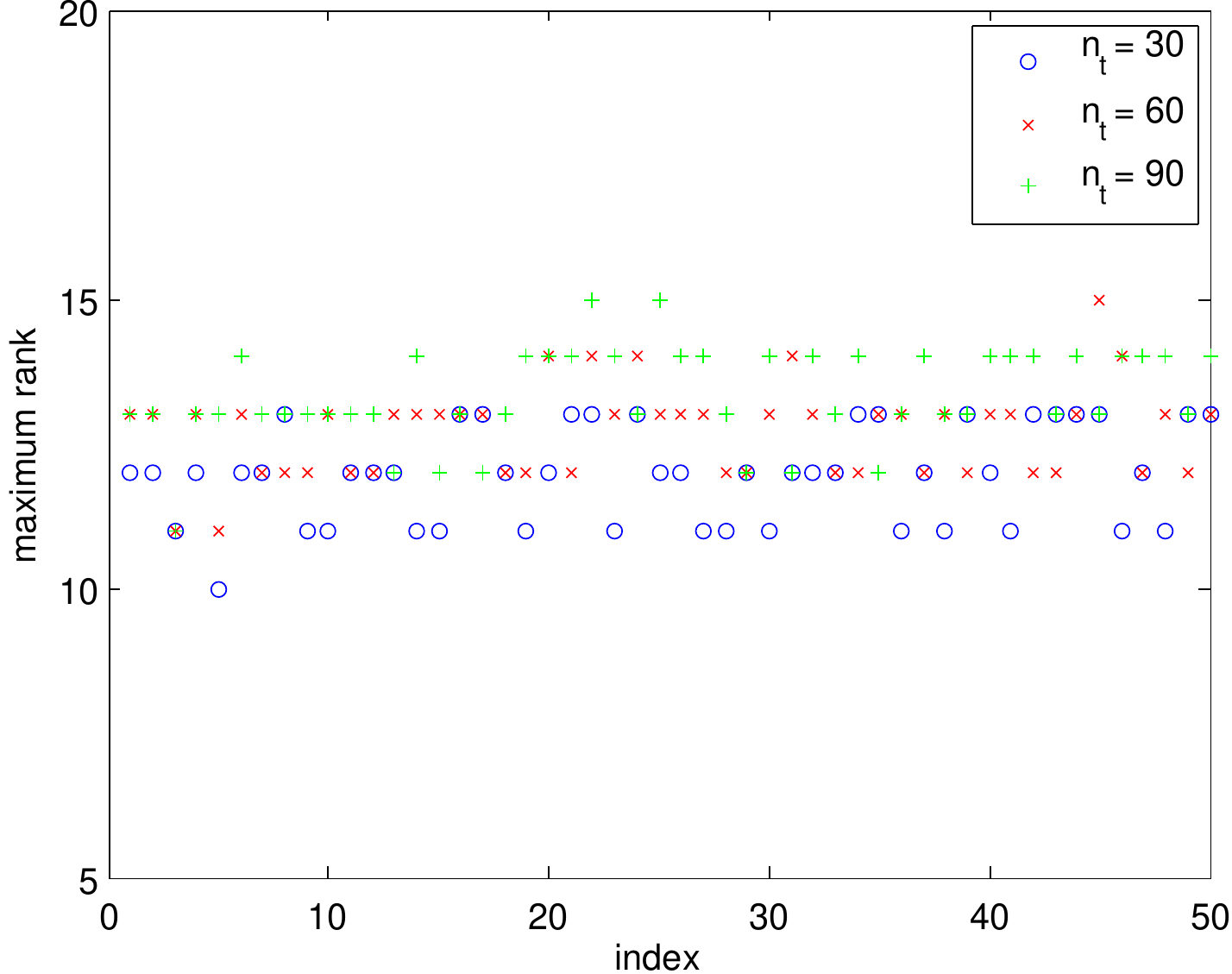}
    \caption{Maximum rank for eigenvectors of $\tilde{\mathcal{H}}_{\text{mis}}.$}\label{fig::rk_l_6_nt369}
\end{figure}

We have seen from~\cref{fig::heat_eig_nt369} that the eigenvalues of $\tilde{\mathcal{H}}_{\text{mis}}$ have a sufficiently fast decay. Meanwhile, each eigenvector also exhibits a low-rank property. Here we plot the maximum rank used to compute the low-rank approximation of eigenvectors corresponding to the $50$ largest eigenvalues for $n_t=30,\ 60,\ 90$, respectively. In~\cref{fig::rk_l_6_nt369}, it is shown that the increase of time steps keeps the rank bounded for the low-rank approximation of the 
eigenvectors of $\tilde{\mathcal{H}}_{\text{mis}}$. The threshold for the low-rank approximation is set to be $10^{-8}$. 

As illustrated in~\Cref{sec::eigenfun}, the eigenvector also admits a low-rank property. We perform a low-rank approximation on each Arnoldi vector throughout the Arnoldi iteration. Since the low-rank approximation is employed, the orthogonality of the basis of Arnoldi vectors is lost.  We just need to compute a few Arnoldi vectors in practice. Here, we use a modified Gram-Schmidt method to perform the full reorthogonalization. Other types of reorthogonalization such as selective orthogonalization~\cite{book::golubvanloan} or periodic orthogonalization~\cite{FlaWAHBWG11} are also possible.

\begin{figure}[htb!]
    \centering
    \subfloat[$n_t=30$]{\includegraphics[width=0.4\textwidth]{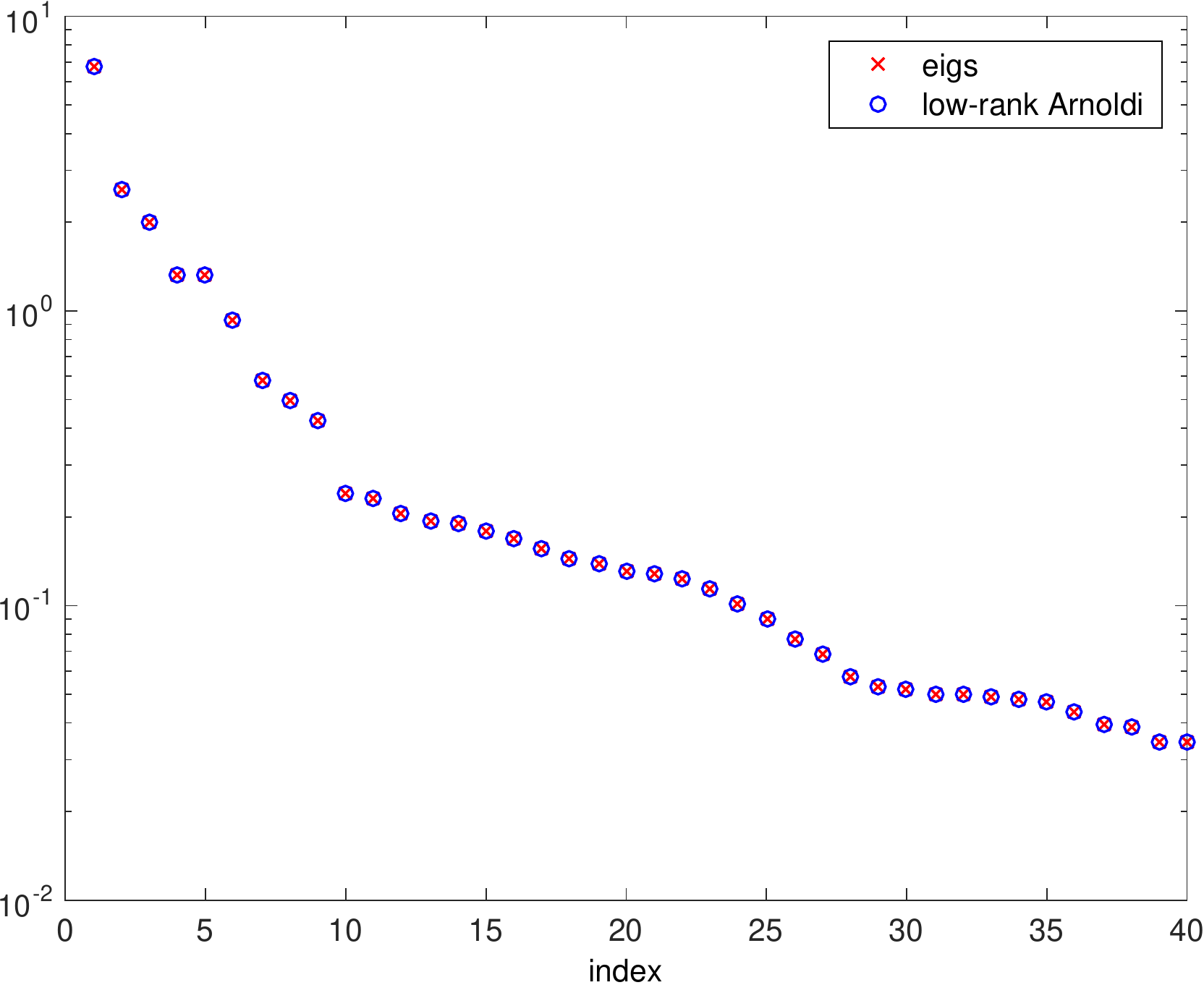}}\qquad
    \subfloat[$n_t=60$]{\includegraphics[width=0.4\textwidth]{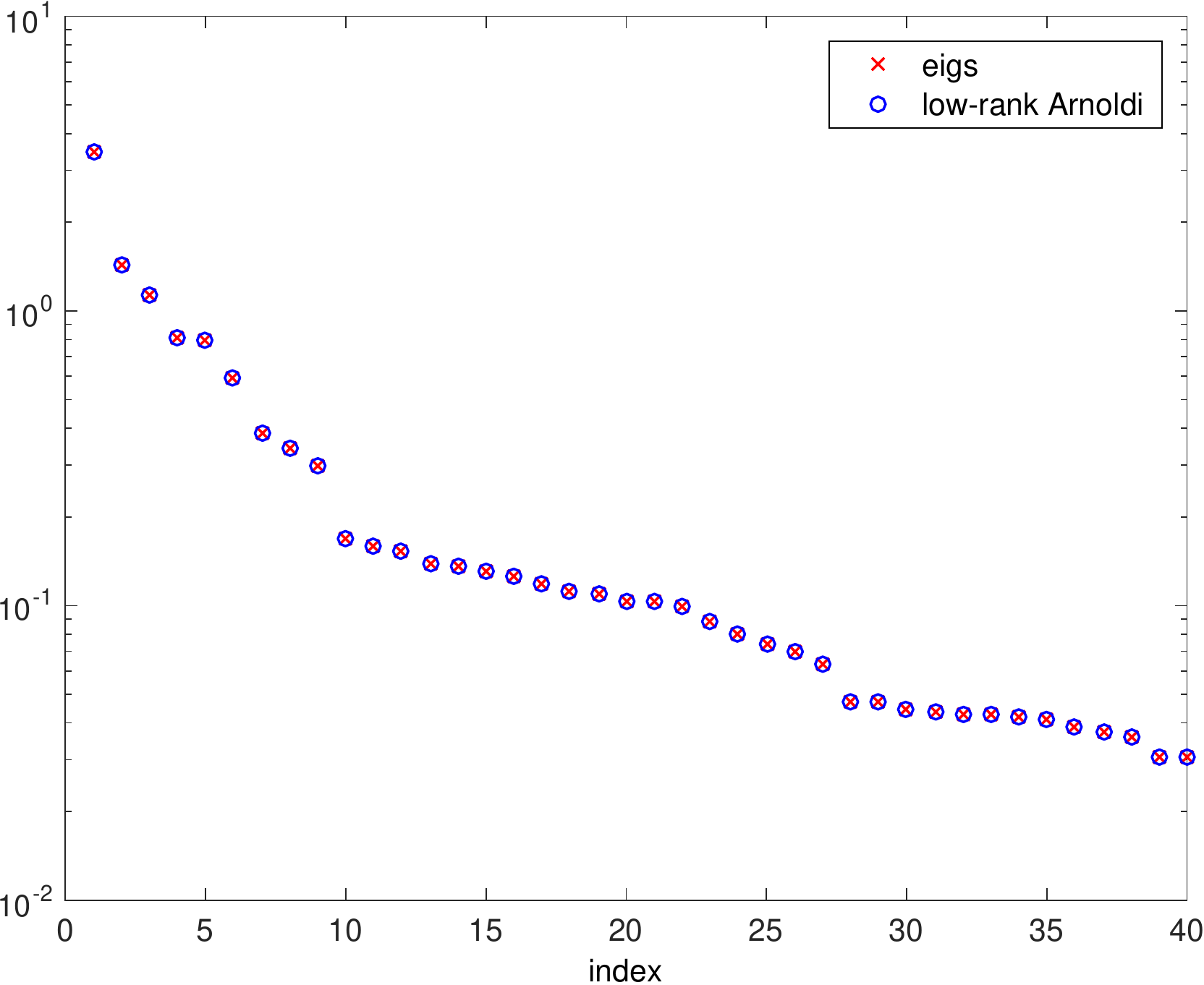}}
    \caption{Eigenvalues of $\tilde{\mathcal{H}}_{\text{mis}}$ computed using `{\texttt eigs}' and low-rank Arnoldi.}\label{fig::heat_eig_lanczos}
\end{figure}

Here we use examples discretized by a $32\times 32$ grid, with $30$ and $60$ time steps to illustrate the effectiveness of the low-rank Arnoldi method. First, we plot the largest $40$ eigenvalues of $\tilde{\mathcal{H}}_{\text{mis}}$ computed using the low-rank Arnoldi method. Next, we compute $\tilde{\mathcal{H}}_{\text{mis}}$ explicitly and use `{\texttt eigs}' in MATLAB to compute its $40$ largest eigenvalues. These results are shown in~\cref{fig::heat_eig_lanczos}. It is clearly shown that the low-rank Arnoldi method can recover the eigenvalues exactly by using reorthogonalization. 

As shown in~\cref{fig::heat_eig_nt369}, the increase of time steps does not influence the decay rate of the eigenvectors of $\tilde{\mathcal{H}}_{\text{mis}}$. Next we will show that an increase of the spatial parameters does not influence the decay rate of eigenvalues of $\tilde{\mathcal{H}}_{\text{mis}}$ either. 

\begin{figure}[H]
    \centering
    \subfloat[largest $60$ eigenvalues of $\tilde{\mathcal{H}}_{\text{mis}}$]{\label{fig::eig_l_69_nt_30}\includegraphics[width=0.4\textwidth]{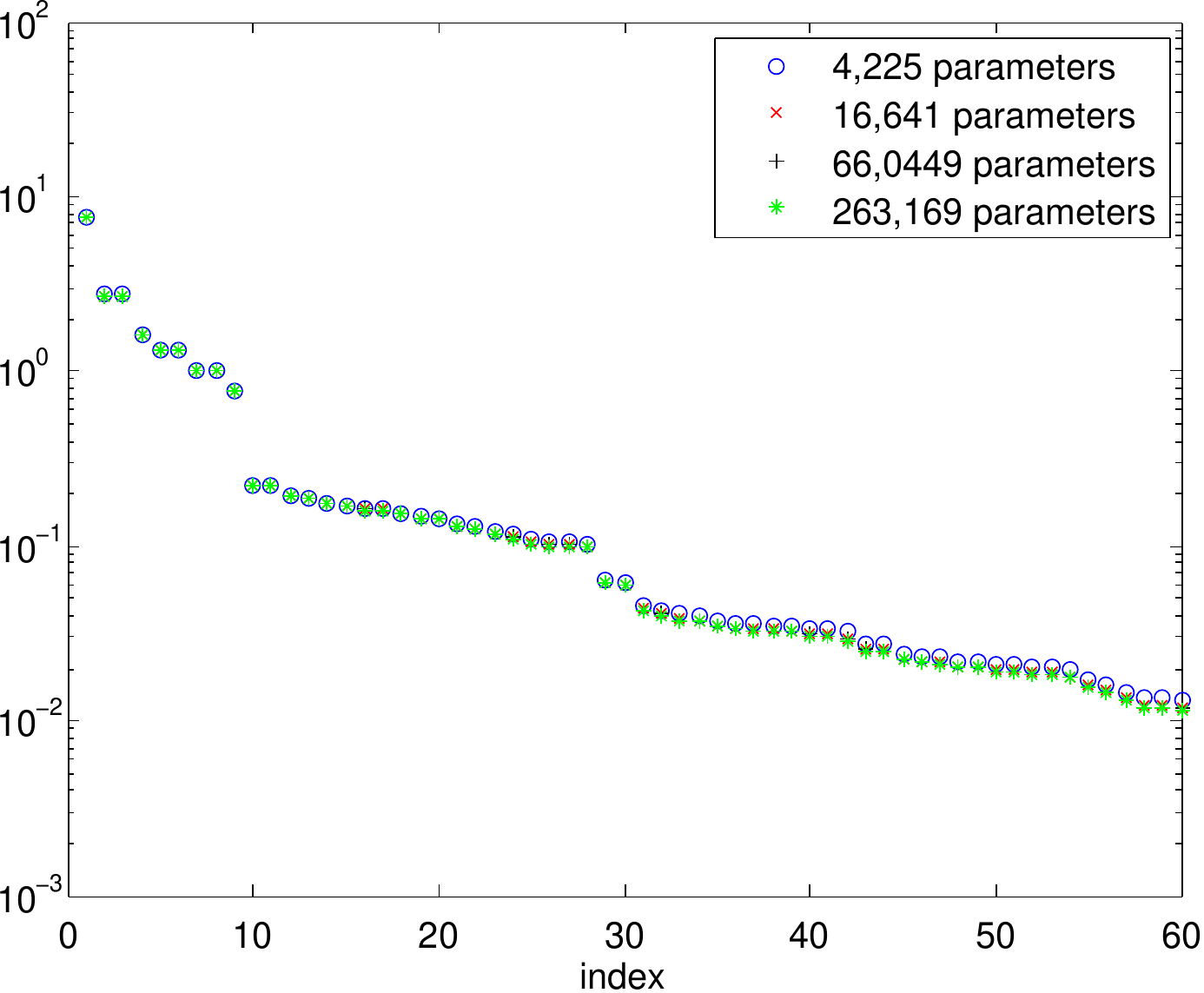}}\qquad
    \subfloat[maximum rank]{\label{fig::max_rk_l_69_nt_30}\includegraphics[width=0.4\textwidth]{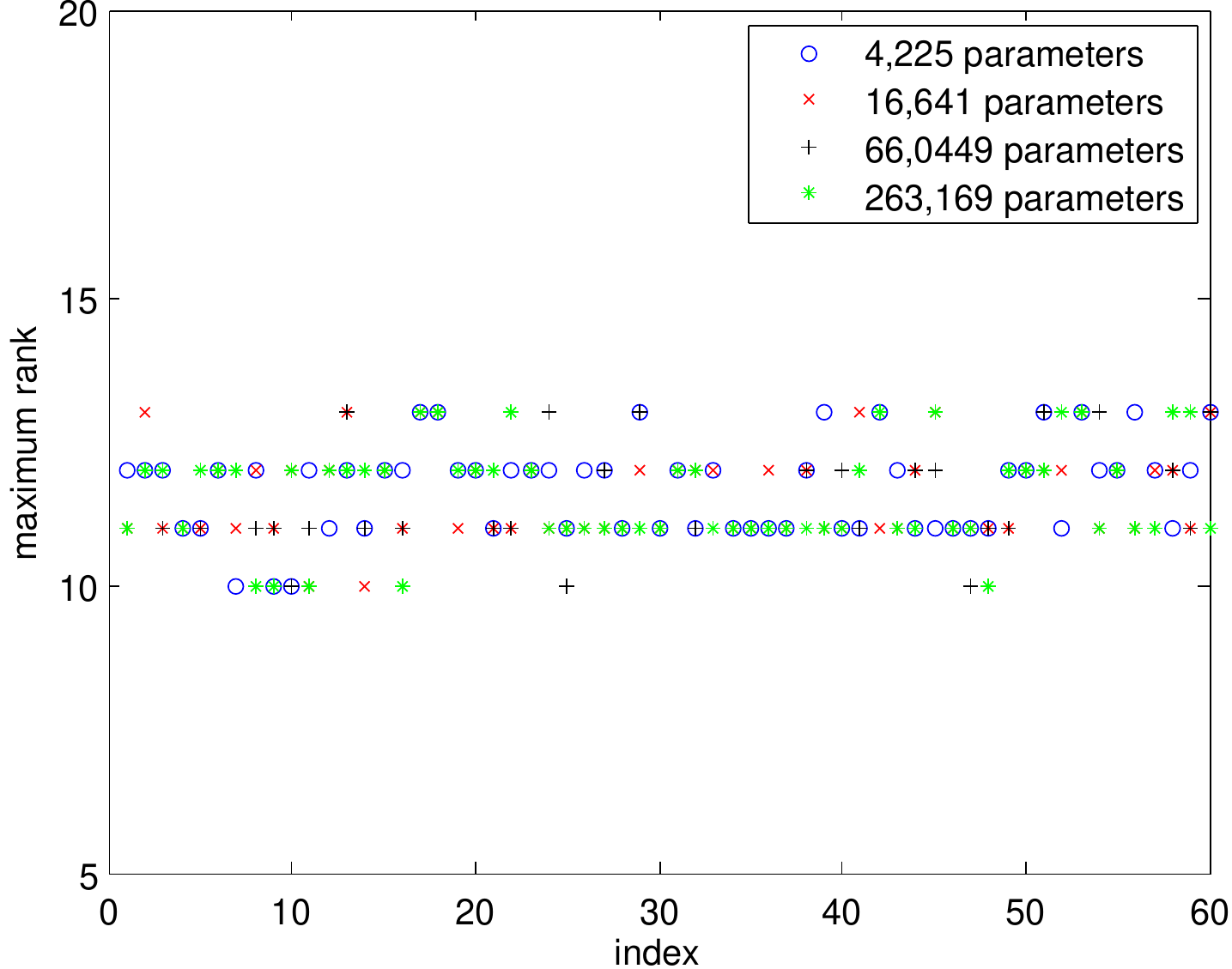}}
    \caption{Eigenvalues for $\tilde{\mathcal{H}}_{\text{mis}}$ and maximum rank for low-rank Arnoldi method, $n_t=30$.}\label{fig::l_69_nt_30}
\end{figure}

We fix the number of time steps $n_t$ to $30$, and compute $60$ largest eigenvalues of $\tilde{\mathcal{H}}_{\text{mis}}$, the results are shown in~\cref{fig::eig_l_69_nt_30}. The maximum rank used for the low-rank Arnoldi method with different number of spatial parameters is shown in~\cref{fig::max_rk_l_69_nt_30}. 

\Cref{fig::eig_l_69_nt_30} shows that the increase number of parameters do not influence the decay property of eigenvalues of $\tilde{\mathcal{H}}_{\text{mis}}$. It is expected that for different number of parameters, the eigenvalues of $\tilde{\mathcal{H}}_{\text{mis}}$ converge to the eigenvalues of the prior-preconditioned operator as long as the discretization of parameter field is good enough. This is illustrated by the eigenvalues shown in~\cref{fig::eig_l_69_nt_30}. Meanwhile,
maximum ranks used in the low-rank Arnoldi method are also bounded by a constant with the increase of number of parameters, which is shown in \cref{fig::max_rk_l_69_nt_30}.

As stated before, a threshold of $10^{-1}$ for the eigenvalue computations of $\tilde{\mathcal{H}}_{\text{mis}}$ is enough to reduce the uncertainty and to approximate the posterior covariance matrix. Next, we plot the diagonal entries of the approximated posterior covariance matrix, i.e., the variance of the points for a $64\times 64$ mesh with a different truncation threshold $\epsilon$ for eigenvalue computations of $\tilde{\mathcal{H}}_{\text{mis}}$. We use a $9$ sensors setting for the sparse observation inverse problem, where $9$ sensors are uniformly distributed inside the domain and the size of each sensor is $\sfrac{1}{256}$ of the domain. Here we set $\beta_{\text{noise}}=10^4\beta_{\text{prior}}$ and the prior covariance matrix $\Gamma_{\text{prior}}=10 I$, where $I$ is an identity matrix with appropriate size.

\begin{figure}[H]
    \centering
    \subfloat[$\epsilon=10^{0}$]{\label{fig::post_variance_a}\includegraphics[width=0.45\textwidth]{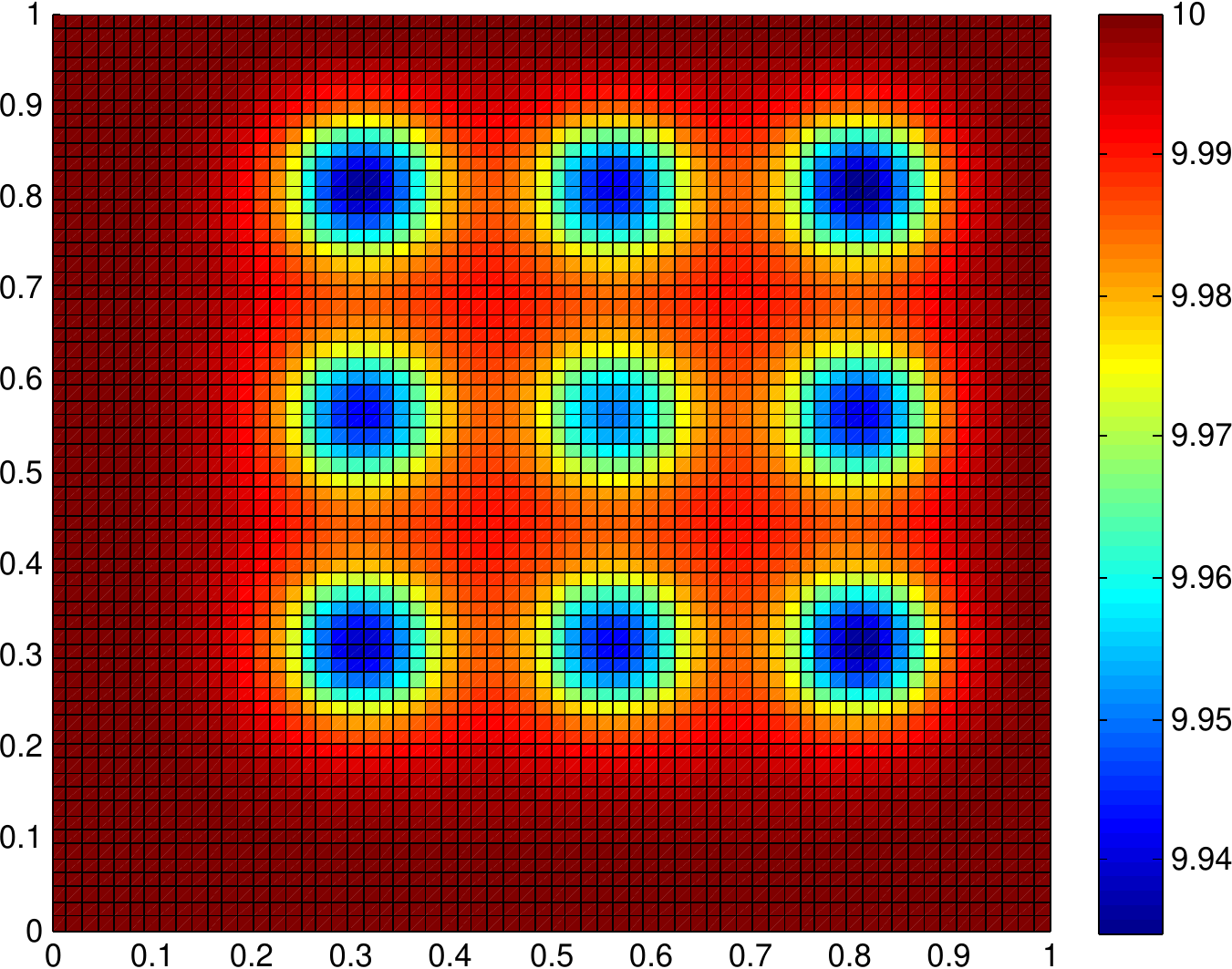}}\qquad
    \subfloat[$\epsilon=10^{-1}$]{\label{fig::post_variance_b}\includegraphics[width=0.45\textwidth]{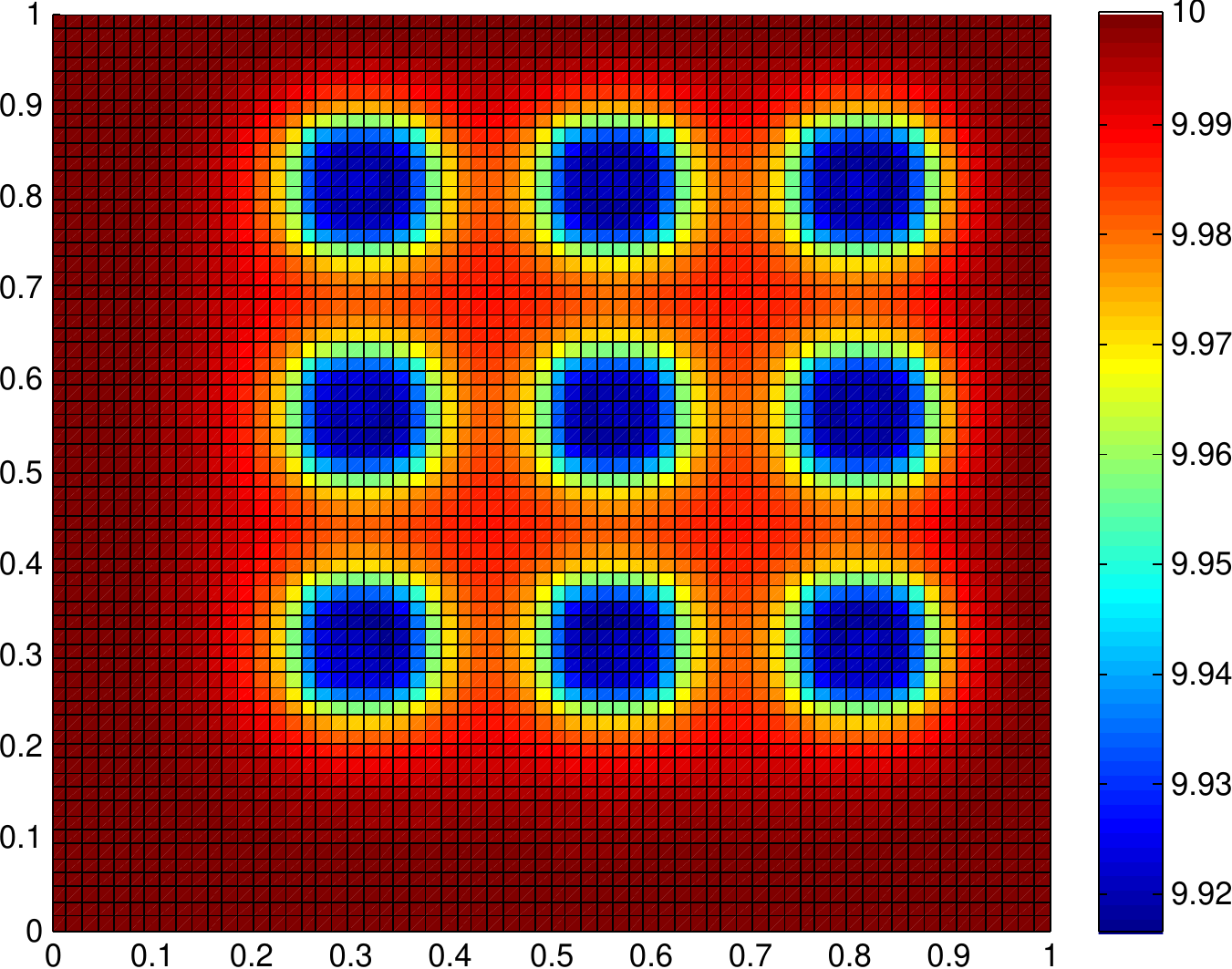}}\\
    \subfloat[$\epsilon=10^{-2}$]{\label{fig::post_variance_c}\includegraphics[width=0.45\textwidth]{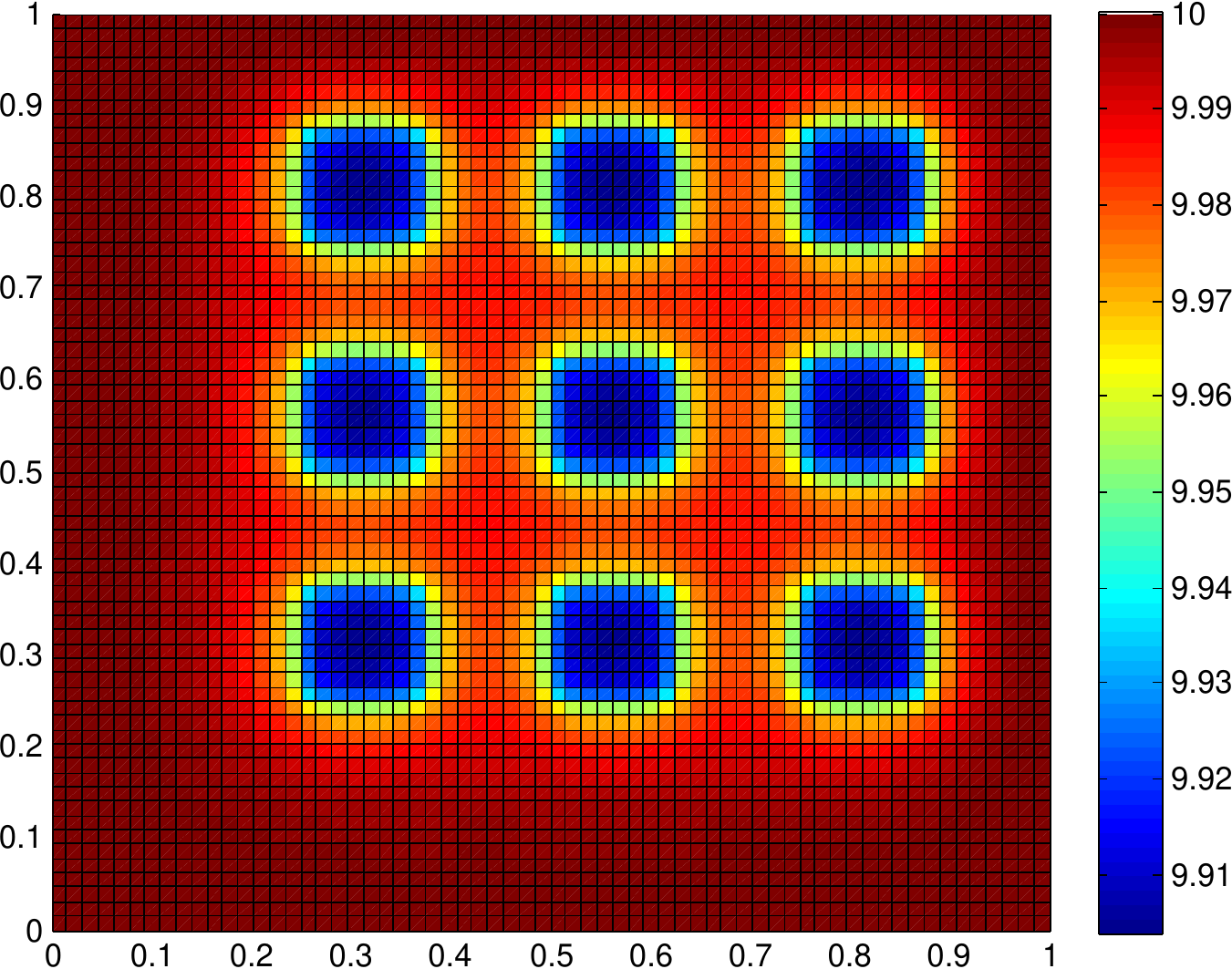}}\qquad
    \subfloat[$\epsilon=10^{-3}$]{\label{fig::post_variance_d}\includegraphics[width=0.45\textwidth]{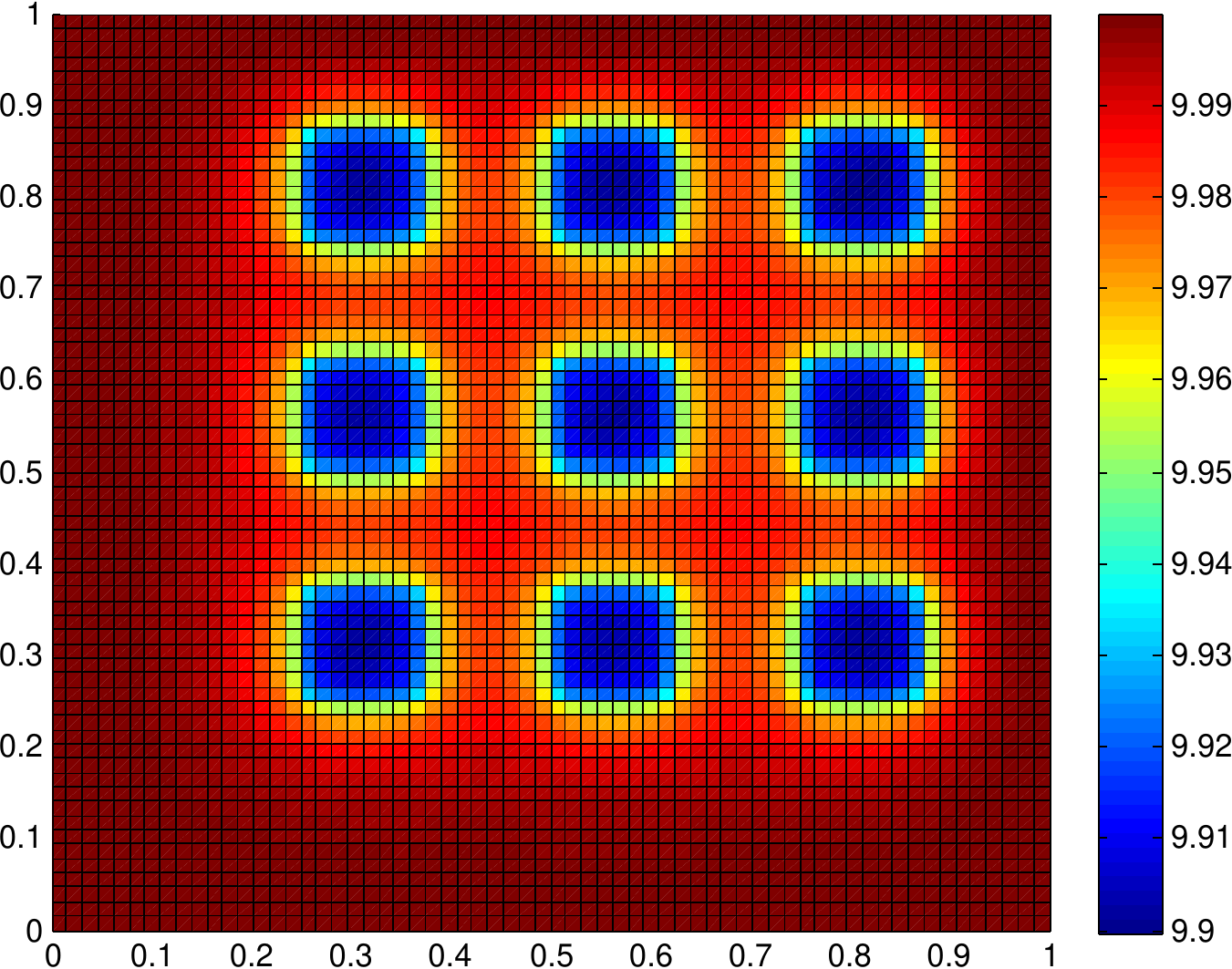}}
    \caption{Diagonal entries of $\Gamma_{\text{post}}$, $n_t=30$, $\beta_{\text{noise}}=10^4\beta_{\text{prior}}.$}\label{fig::post_variance}
\end{figure}

For the sparse observation case, covariance updates are mostly clustered around the area where data are observed while the rest are dominated by the prior.  Uncertainty can only be reduced at areas around the location of sensors. This is clearly shown by~\cref{fig::post_variance_a}~-~\cref{fig::post_variance_d}, where the dark colored areas are placed at the location of the sensors and have the lowest variance. Decreasing the threshold $\epsilon$, we observe that the variance is further reduced around the location of sensors. For smaller values of $\epsilon$ no more reduction in the variance is achieved as all essential information are already captured. \Cref{fig::post_variance} shows that as long as the computations of $\Gamma_{\text{post}}$ is convergent, using a threshold $\epsilon=10^{-1}$ is enough to approximate the posterior covariance matrix and to reduce the uncertainty.

\begin{figure}[H]
    \centering
    \subfloat[$\epsilon=10^{0}$]{\label{fig::post_variance_a_6}\includegraphics[width=0.45\textwidth]{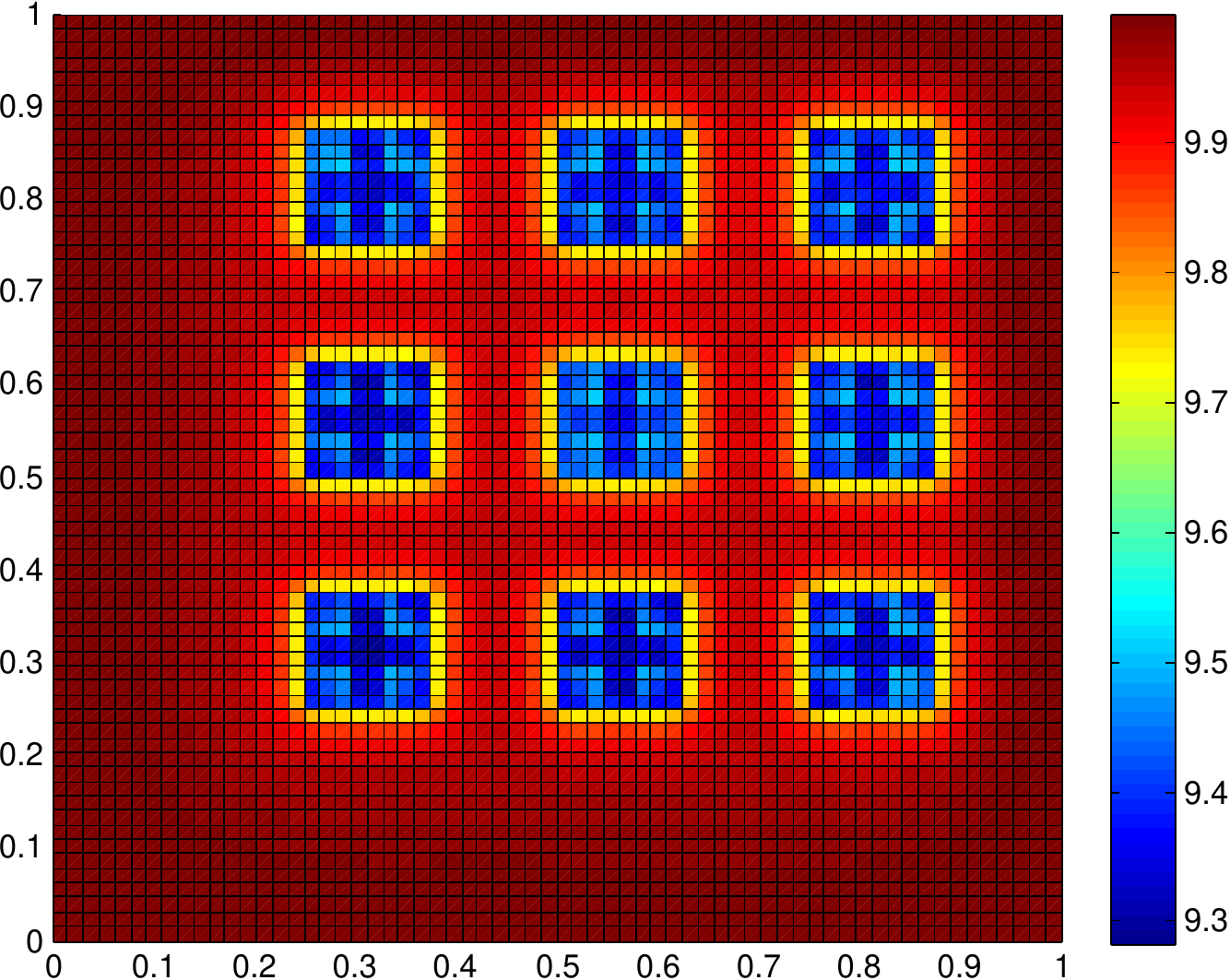}}\quad
    \subfloat[$\epsilon=10^{-1}$]{\label{fig::post_variance_b_6}\includegraphics[width=0.45\textwidth]{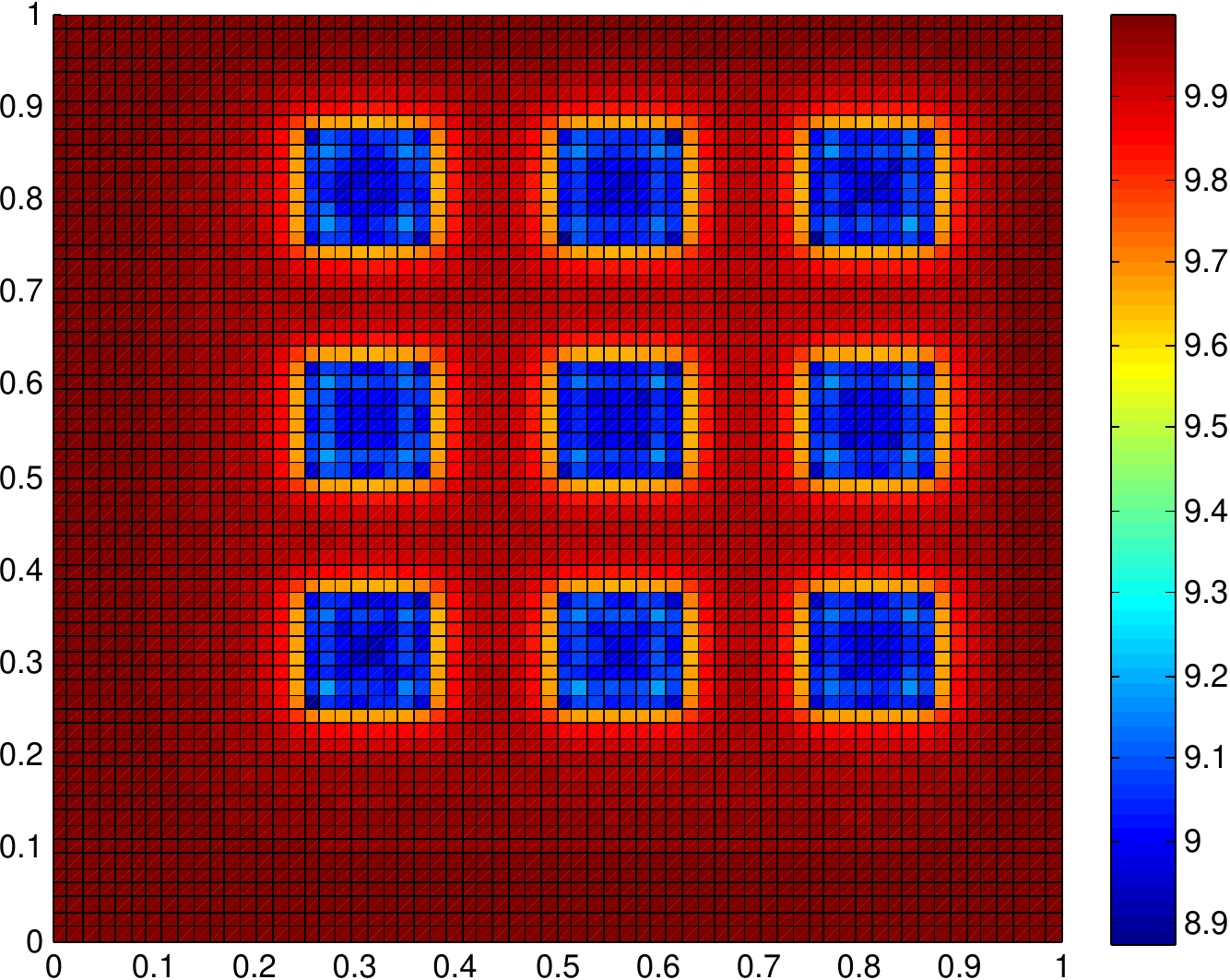}}
    \caption{Diagonal entries of $\Gamma_{\text{post}}$, $n_t=30$, $\beta_{\text{noise}}=10^6\beta_{\text{prior}}$.}\label{fig::post_variance_6}
\end{figure}

\begin{figure}[H]
    \centering
    \subfloat[$\epsilon=10^{0}$]{\label{fig::post_variance_a_8}\includegraphics[width=0.45\textwidth]{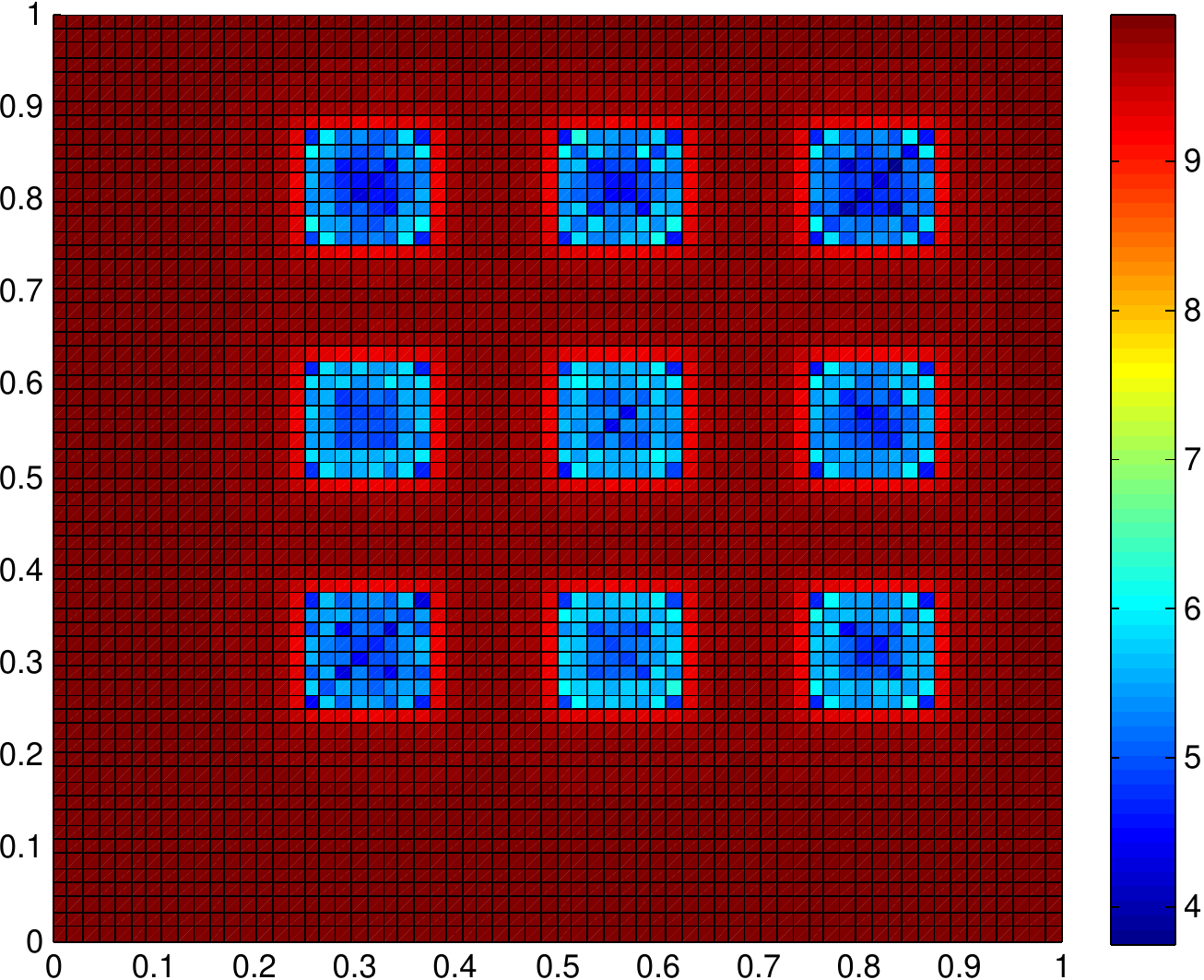}}\quad
    \subfloat[$\epsilon=10^{-1}$]{\label{fig::post_variance_b_8}\includegraphics[width=0.45\textwidth]{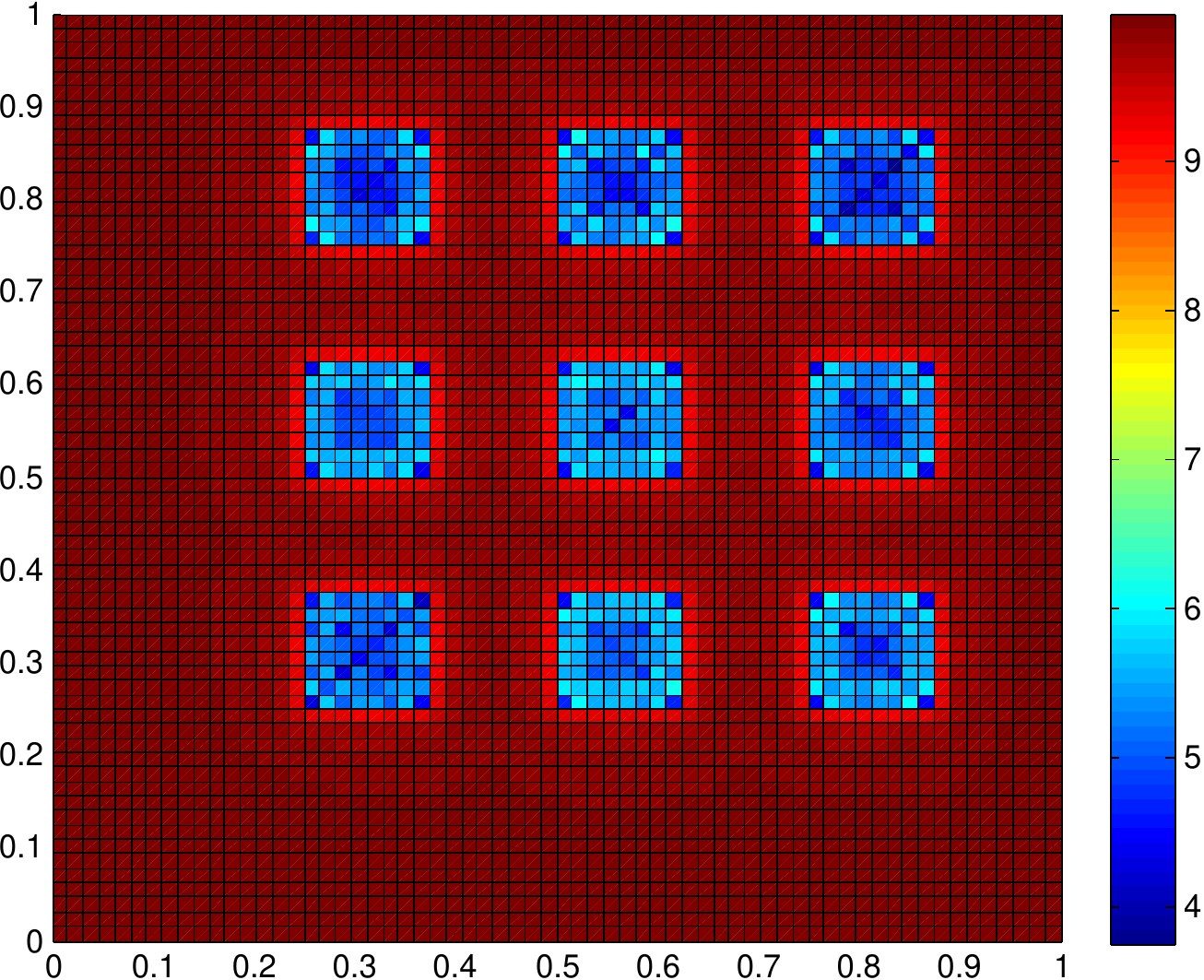}}
    \caption{Diagonal entries of $\Gamma_{\text{post}}$, $n_t=30$, $\beta_{\text{noise}}=10^8\beta_{\text{prior}}$}\label{fig::post_variance_8}
\end{figure}

As analyzed in~\Cref{sec::eigenfun}, the eigenvalues of $\tilde{\mathcal{H}}_{\text{mis}}$ are related to $\frac{\beta_{\text{noise}}}{\beta_{\text{prior}}}$ and this ratio gives different updates of the posterior covariance matrix. Next, we use different $\frac{\beta_{\text{noise}}}{\beta_{\text{prior}}}$ ratios to plot the variance of the parameters. These results are shown in~\Cref{fig::post_variance_6}~-~\Cref{fig::post_variance_8}. The prior covariance matrix is set to be $\Gamma_{\text{prior}}=10 I$, where $I$ is an identity matrix with appropriate size.

For the case $\beta_{\text{noise}}=10^6\beta_{\text{prior}}$, we need $72$ Arnoldi iterations for $\epsilon=10^{0}$ while $163$ Arnoldi iterations for $\epsilon=10^{-1}$. With $\beta_{\text{noise}}=10^8\beta_{\text{prior}}$, we need $347$ Arnoldi iterations for $\epsilon=10^{0}$ while $470$ Arnoldi iterations for $\epsilon=10^{-1}$.

\cref{fig::post_variance_6} and \cref{fig::post_variance_8} show that with the increase of the ratio between $\beta_{\text{noise}}$ and $\beta_{\text{prior}}$, the diagonal entries of $\Gamma_{\text{post}}$ becomes smaller. This implies that to further reduce the posterior variance, we need bigger ratio between $\beta_{\text{noise}}$ and $\beta_{\text{prior}}$. This can be explained as follows, the weight for the data misfit part in the optimization problem~\eqref{eqn::obj} is getting bigger for bigger  $\frac{\beta_{\text{noise}}}{\beta_{\text{prior}}}$. This means that the data misfit part is more strictly optimized than for smaller $\frac{\beta_{\text{noise}}}{\beta_{\text{prior}}}$. Therefore the error between the estimation and observed data is getting smaller. However, this yields a more ill-conditioned problem and more Arnoldi iterations are needed. Therefore, a balance between covariance reduction and computational effort is needed with our approach enabling the storage of many Arnoldi vectors due to 
the complexity reduction of the low-rank approach.

\subsection{Convection-Diffusion Equation}

In this section, we study our low-rank approach for a stochastic inverse convection-diffusion problem. Here, the convection-diffusion operator $\mathcal{L}$ is given by

\[
    \mathcal{L}= -\nu \Delta u + \overrightarrow{\omega}\cdot \nabla u. 
\]

The computational domain is chosen as a square domain given by $[0,\ 1]\times [0,\ 1]$, $\overrightarrow{\omega}=(0,\ 1)$, and the inflow is posed on the down boundary while the outflow is posed on the upper boundary. Boundary conditions are prescribed according to the analytic solution of the convection-diffusion equation, which is described as Example 3.3.1 in~\cite{book::andy}. We use the streamline upwind Petrov–Galerkin (SUPG) finite element method to discretize the convection-diffusion equation.

\begin{figure}[H]
	\centering
	\subfloat[$50$ largest eigenvalues of $\tilde{\mathcal{H}}_{\text{mis}}$]{\label{fig::cd_eig_lanczos}\includegraphics[width=0.46\textwidth]{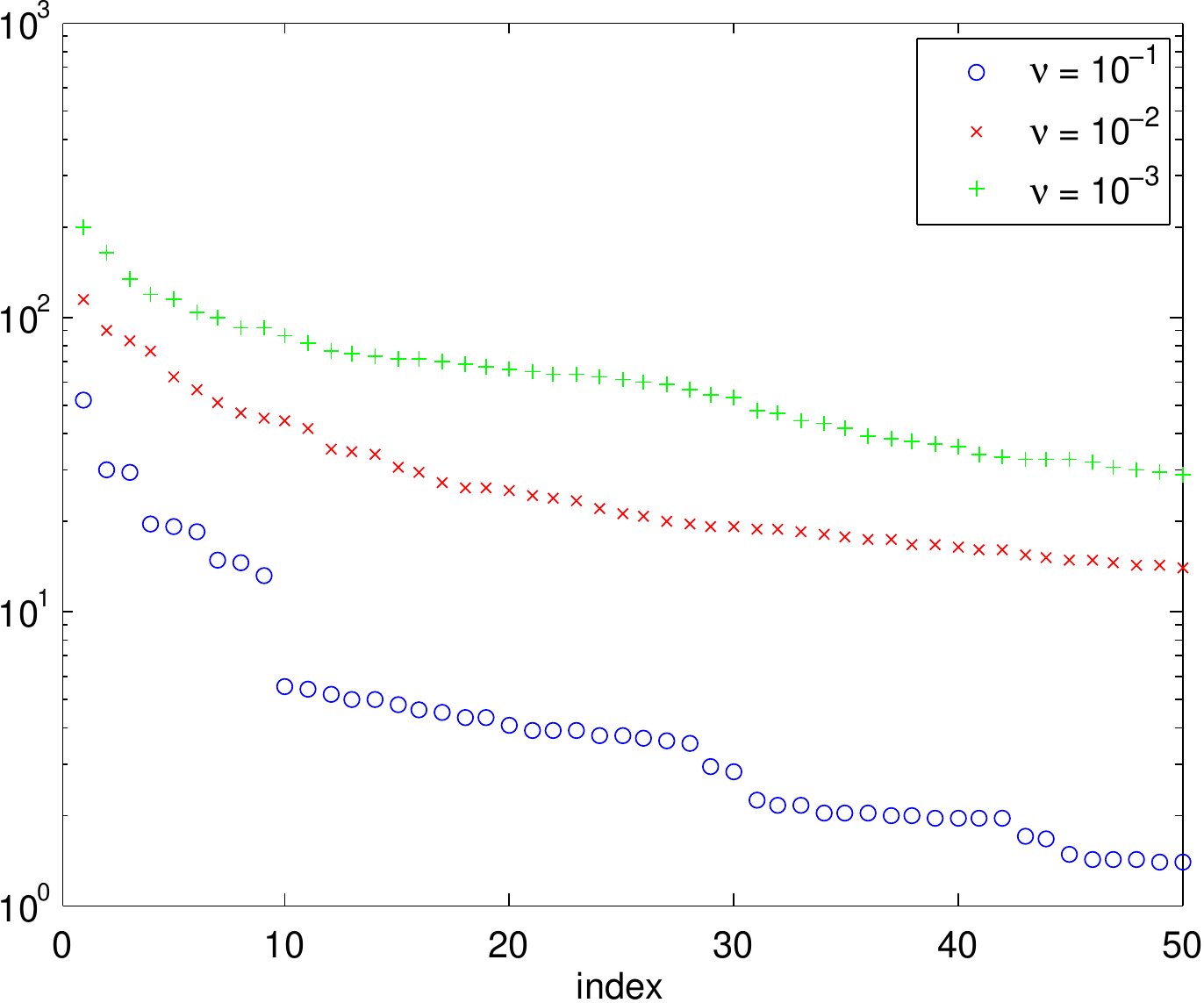}}\quad
	\subfloat[maximum ranks at each low-rank Arnoldi iteration]{\label{fig::cd_rk_lanczos}\includegraphics[width=0.48\textwidth]{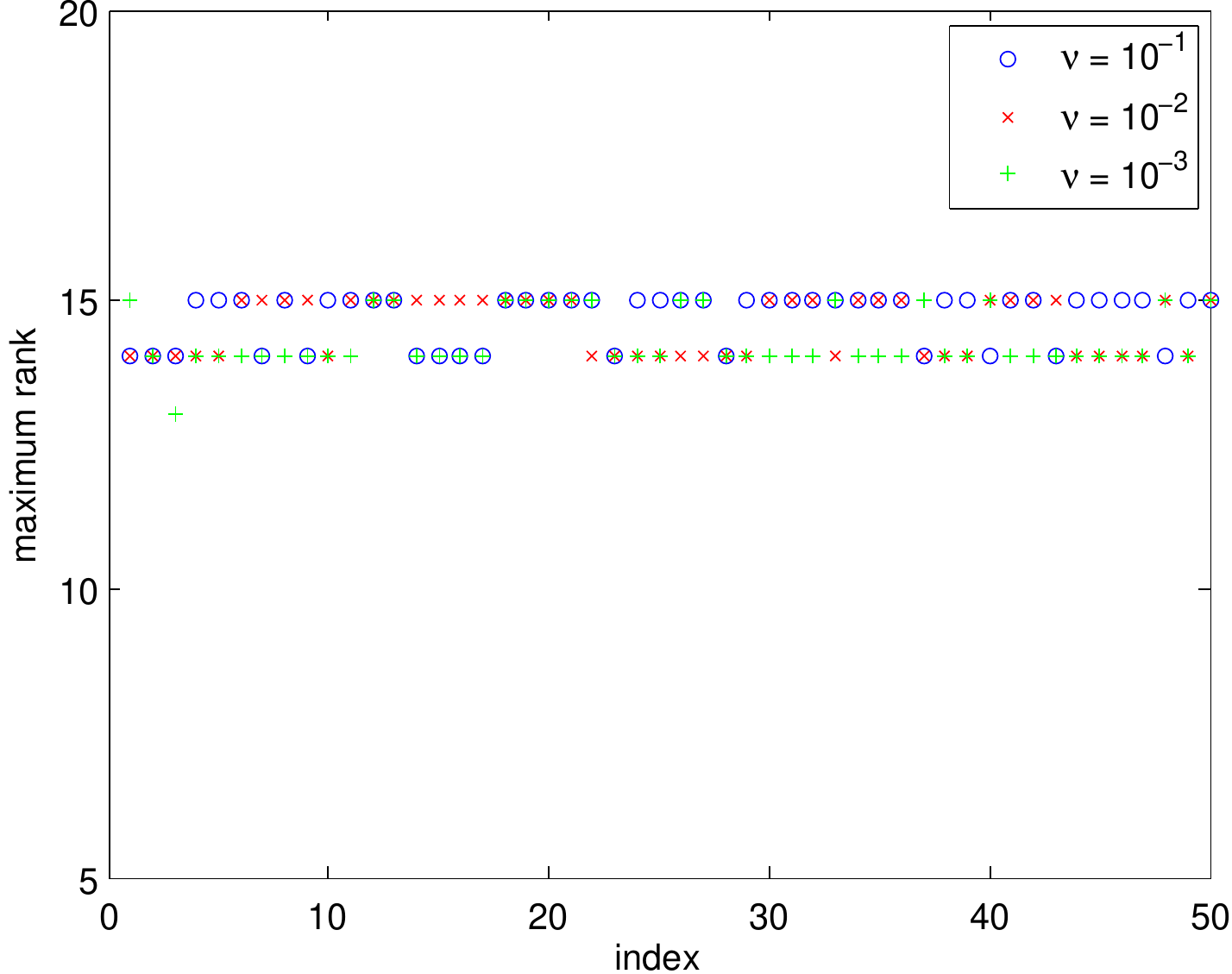}}
	\caption{$50$ largest eigenvalues of $\tilde{\mathcal{H}}_{\text{mis}}$ and maximum ranks at each low-rank Arnoldi iteration with different $\nu$.}\label{fig::cd_e_123_nt_3}
\end{figure}

First, we show the eigenvalue decay of $\tilde{\mathcal{H}}_{\text{mis}}$ for different settings of the viscosity parameter $\nu$. Here we set the number of time steps $n_t$ to be 30, and $\beta_{\text{noise}}=10^{4}\beta_{\text{prior}}$. We plot the $50$ largest eigenvalues of $\tilde{\mathcal{H}}_{\text{mis}}$ for different $\nu$ in~\cref{fig::cd_eig_lanczos}.

As shown by~\cref{fig::cd_eig_lanczos}, the eigenvalues of $\tilde{\mathcal{H}}_{\text{mis}}$ decay rapidly for big $\nu$ while this decay rate slows down when $\nu$ gets smaller. Therefore, more Arnoldi iterations are needed to get a satisfactory approximation of $\tilde{\mathcal{H}}_{\text{mis}}$. For smaller $\nu$, the largest eigenvalue is also bigger than that for bigger $\nu$ as shown in~\cref{fig::cd_eig_lanczos}. The first few eigenvectors form a more dominant subspace than for bigger $\nu$. It is therefore possible to choose a larger truncation threshold for smaller $\nu$, which will be shown later.


The maximum rank for the low-rank Arnoldi iteration with different $\nu$ is shown in~\cref{fig::cd_rk_lanczos}. It shows that the maximum rank does not increase with the decrease of $\nu$ and is bounded by a small constant. Therefore, the complexity for both computations and storage is $\mathcal{O}(n_x + n_t)$ for the low-rank Arnoldi approach.  

As shown in~\cref{fig::cd_eig_lanczos}, the eigenvalues of $\tilde{\mathcal{H}}_{\text{mis}}$ decay slower for smaller $\nu$, and more Arnoldi iterations are needed, which is due to the property of the problem. For such a problem with smaller $\nu$, our low-rank approach is much more superior to the standard Arnoldi method introduced in~\cite{FlaWAHBWG11} since we need to compute and store more Arnoldi vectors. Note that this is doable with the approach presented here.    

\begin{figure}[H]
\centering
\subfloat[$50$ largest eigenvalues of $\tilde{\mathcal{H}}_{\text{mis}}$]{\label{fig::cd_eig_lanczos_nt_369}\includegraphics[width=0.46\textwidth]{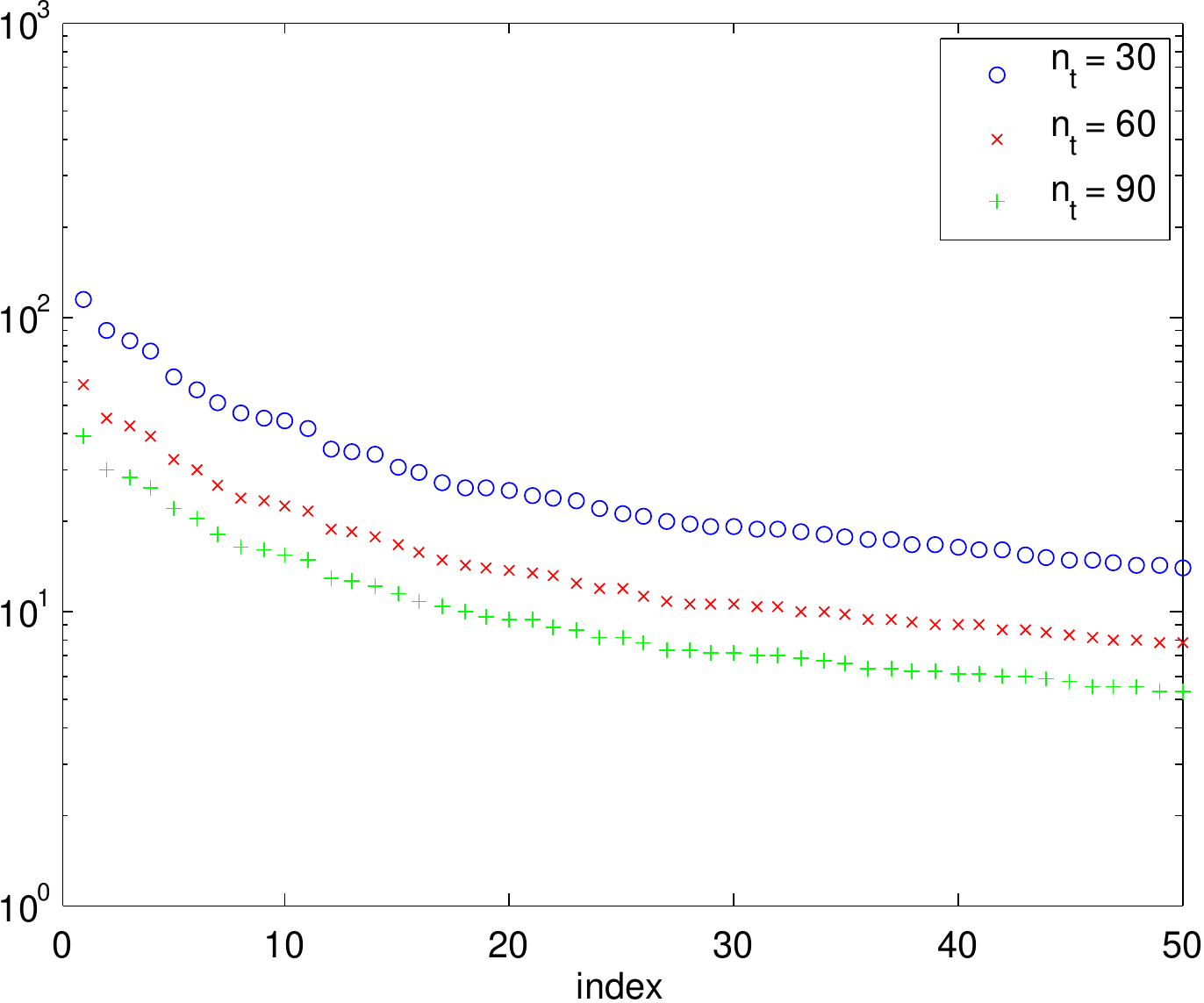}}\quad
\subfloat[maximum ranks at each low-rank Arnoldi iteration]{\label{fig::cd_rk_lanczos_nt_369}\includegraphics[width=0.48\textwidth]{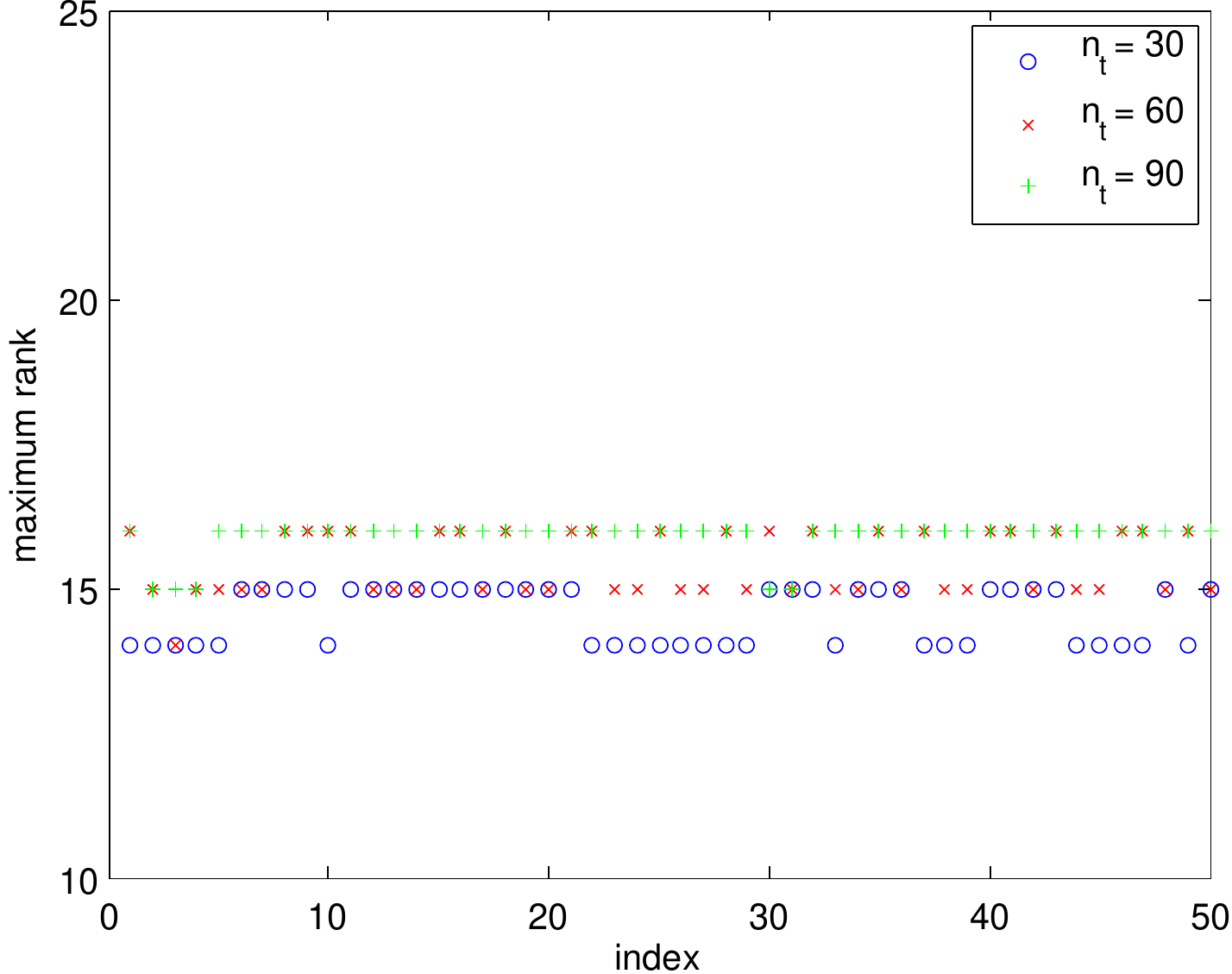}}
\caption{$50$ largest eigenvalues of $\tilde{\mathcal{H}}_{\text{mis}}$ and maximum ranks at each low-rank Arnoldi iteration with different $n_t.$}\label{fig::cd_e_2_nt_369}
\end{figure}

Next, we set $\nu$ to be $10^{-2}$ and compute the $50$ largest eigenvalues of $\tilde{\mathcal{H}}_{\text{mis}}$ with different $n_t$, which are shown in~\cref{fig::cd_eig_lanczos_nt_369}. We are also interested in the relation between the maximum rank at each Arnoldi iteration and the number of time steps $n_t$. This is shown in~\cref{fig::cd_rk_lanczos_nt_369}.

%
%
%

As shown by~\cref{fig::cd_eig_lanczos_nt_369}, with the increase in the number of time steps, the eigenvalue decay of $\tilde{\mathcal{H}}_{\text{mis}}$ behaves similar. Maximum ranks at each low-rank Arnoldi iteration with different $n_t$ are bounded by a moderately small constant, which is independent of $n_t$. 

\cref{fig::cd_rk_lanczos} and \cref{fig::cd_rk_lanczos_nt_369} show that the maximum rank for each low-rank Arnoldi iteration is almost invariant w.r.t.\ the number of time steps $n_t$ and the viscosity parameter $\nu$. This makes our low-rank Arnoldi method quite appealing for even complicated stochastic convection dominated inverse problems over a long time horizon.

\begin{figure}[H]
    \centering
    \subfloat[$\epsilon=10^{1}$]{\label{fig::cd_post_variance_a_8}\includegraphics[width=0.45\textwidth]{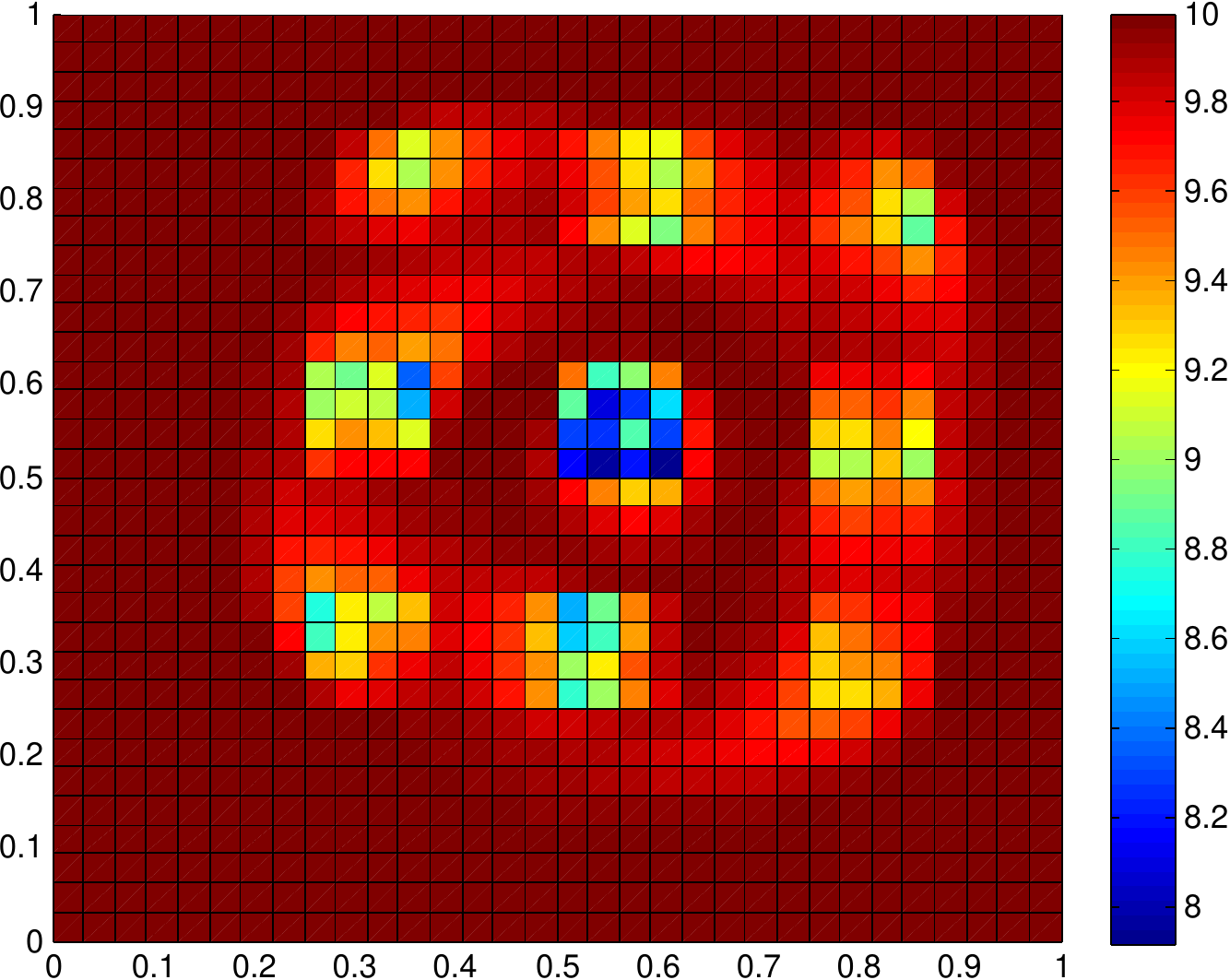}}\quad
    \subfloat[$\epsilon=10^{0}$]{\label{fig::cd_post_variance_b_8}\includegraphics[width=0.45\textwidth]{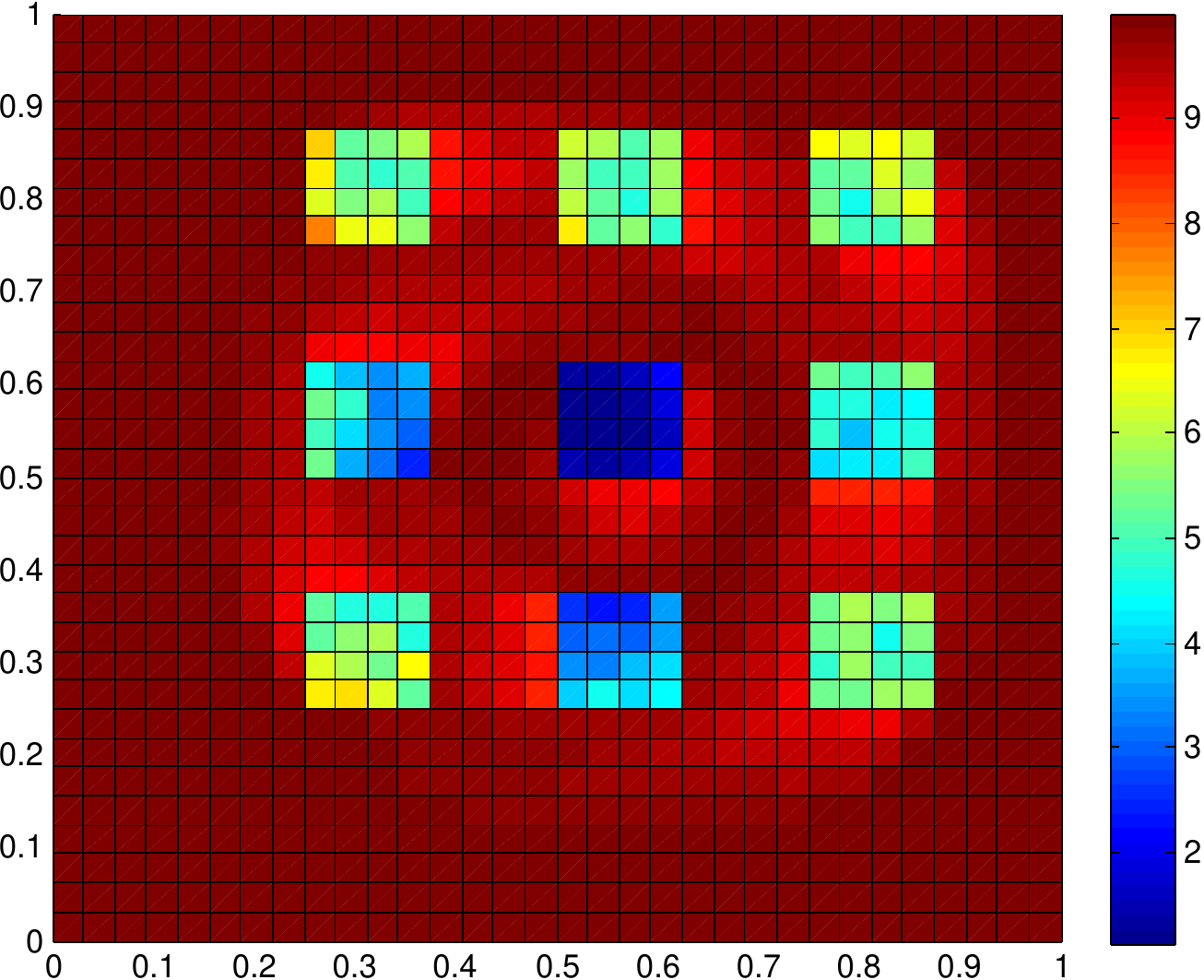}}
    \caption{Diagonal entries of $\Gamma_{\text{post}}$, $\nu=10^{-2}.$}\label{fig::cd_post_variance}
\end{figure}

\begin{figure}[H]
    \centering
    \subfloat[$\epsilon=10^{1}$]{\label{fig::cd_post_variance_a_8_mu_3}\includegraphics[width=0.45\textwidth]{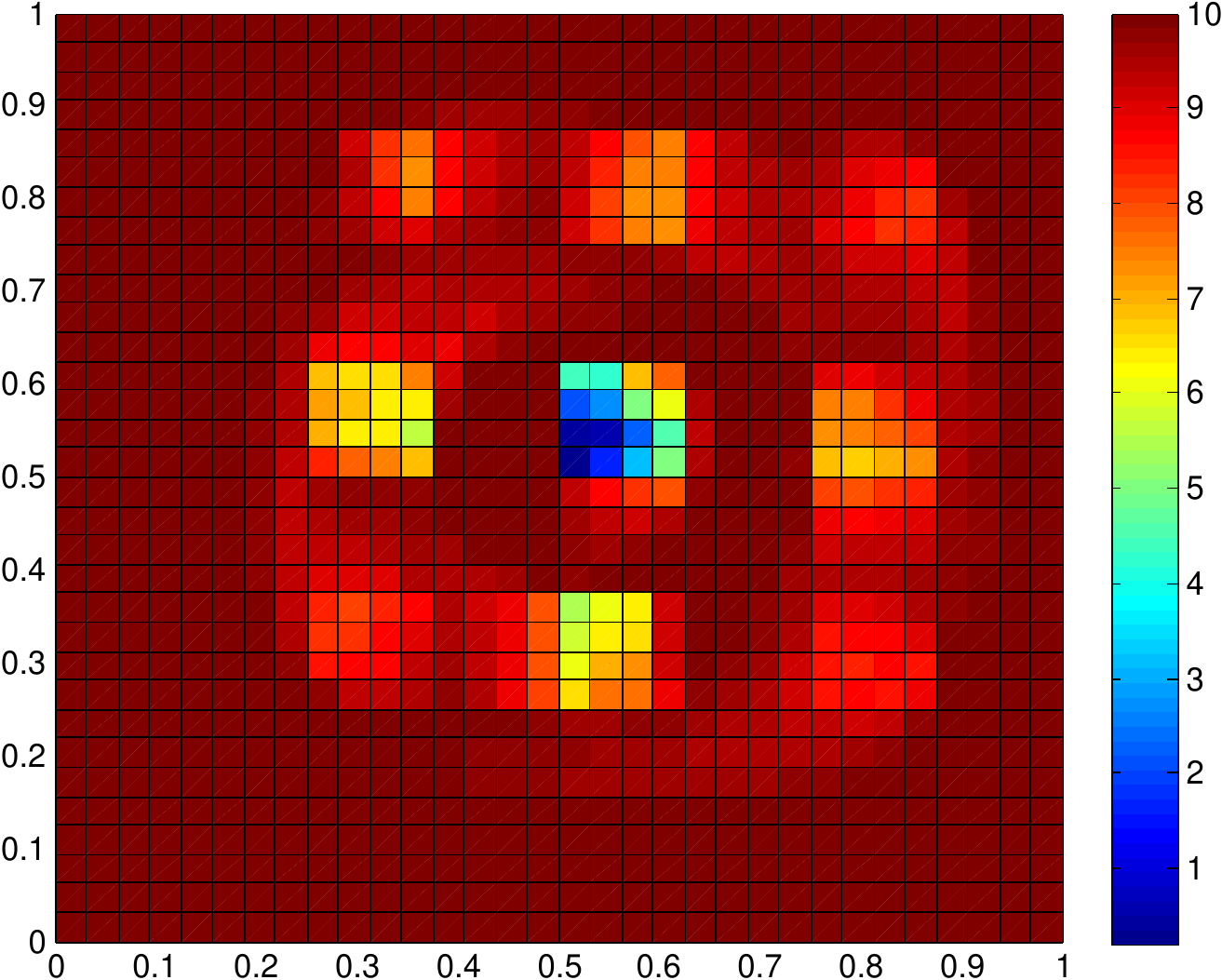}}\quad
    \subfloat[$\epsilon=10^{0}$]{\label{fig::cd_post_variance_b_8_mu_3}\includegraphics[width=0.45\textwidth]{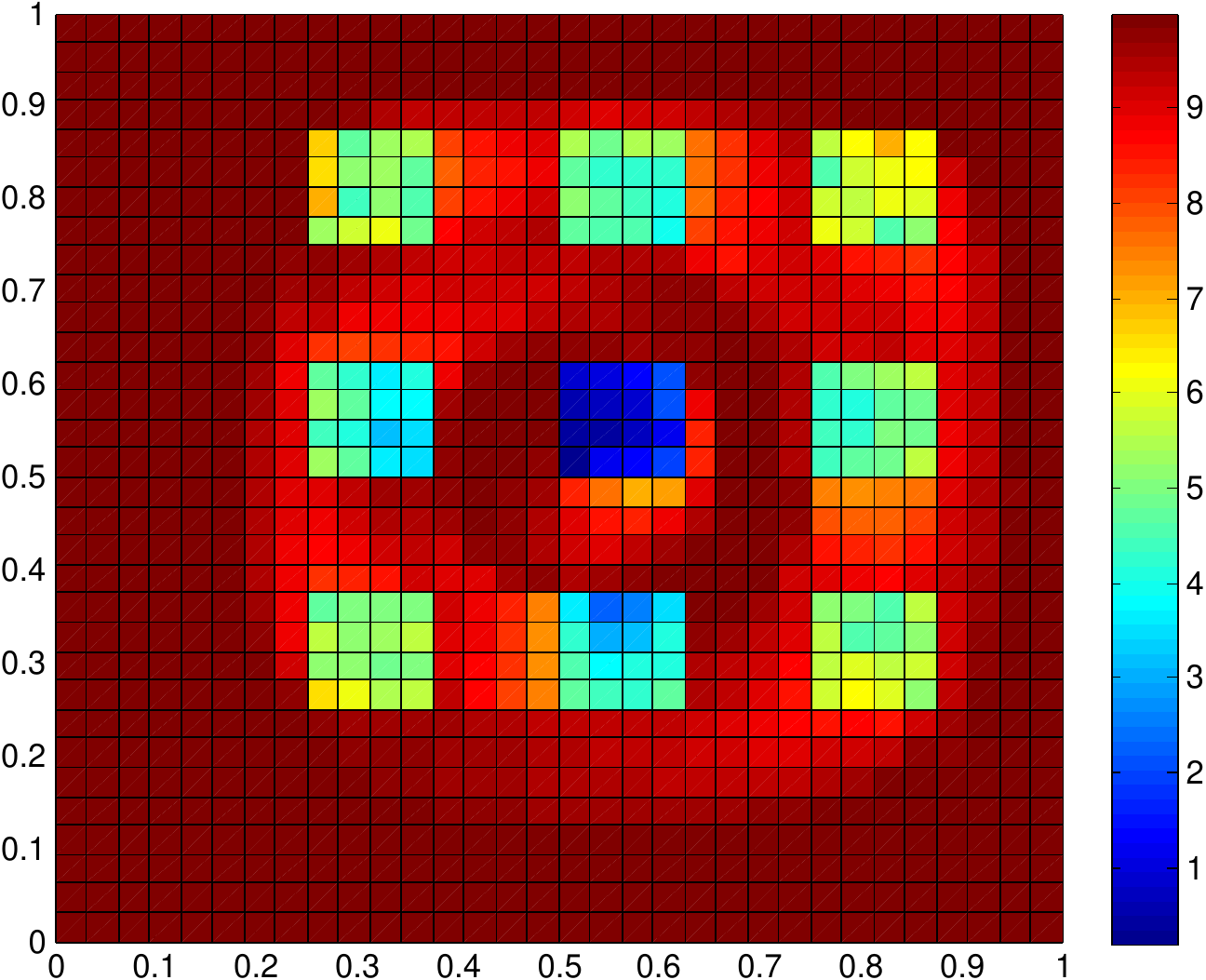}}
    \caption{Diagonal entries of $\Gamma_{\text{post}}$, $\nu=10^{-3}.$}\label{fig::cd_post_variance_mu_3}
\end{figure}

Next, we show the diagonal entries of $\Gamma_{\text{post}}$ for different settings of $\nu$ and the threshold ($\epsilon$) of eigenvalues truncation. Here we set the number of time steps $n_t$ to be $90$, use a $32\times 32$ uniform grid to discretize the convection-diffusion equation, and $\beta_{\text{noise}}=10^{4}\beta_{\text{prior}}$. The results are given by \cref{fig::cd_post_variance} and \cref{fig::cd_post_variance_mu_3}.

For the case $\nu=10^{-2}$, we need $17$ Arnoldi iterations when we use a threshold $\epsilon=10^{1}$, while we need $121$ Arnoldi iterations for $\epsilon=10^{0}$. For the case $\nu=10^{-3}$, we need $52$ low-rank Arnoldi iterations by setting $\epsilon=10^{1}$ and $131$ low-rank Lanczos iterations when we use $\epsilon=10^{0}$. Note that often a further reduction in the truncation parameter does not yield better results as all the essential information is already captured.

\cref{fig::cd_post_variance_b_8} and \cref{fig::cd_post_variance_b_8_mu_3} illustrate that the uncertainty (variance of unknowns) is already reduced dramatically even if we choose a relatively large threshold.

As shown in~\cref{fig::cd_eig_lanczos}, the smaller $\nu$ becomes, the bigger the largest eigenvalue of $\tilde{\mathcal{H}}_{\text{mis}}$ is going to be. Therefore, the first few eigenvectors for smaller $\nu$ form a more dominant subspace. This in turn implies that the uncertainty is much more reduced for smaller $\nu$ when we use the same truncation threshold of eigenvalues. We observe this in \cref{fig::cd_post_variance} and \cref{fig::cd_post_variance_mu_3}.

We can conclude that for the stochastic convection-diffusion inverse problem, our low-rank Arnoldi approach is very flexible and efficient for different time horizon lengths and viscosity parameters. It is even preferred for convection dominated stochastic inverse problems with long time horizon. 

\section{Conclusions}
In this manuscript, we propose a low-rank Arnoldi method to approximate the posterior covariance matrix that appears in stochastic inverse problems. Compared with the standard Arnoldi approach, our approach exploits the low-rank property of each Arnoldi vector and makes a low-rank approximation of such a vector. This reduces the complexity for both computations and storage demand from $\mathcal{O}(n_x n_t)$ to $\mathcal{O}(n_x + n_t)$. Here $n_x$ is the degree of freedom in space and $n_t$ is the degree of freedom in time. This makes solving large scale stochastic inverse problems possible.

Our low-rank approach introduced in this manuscript solves linear stochastic inverse problems that can be put into the Bayesian framework. The next step of our work is to extend the low-rank approach introduced in this manuscript to nonlinear stochastic inverse problems, which is still a big challenge. 
\section*{Acknowledgements}
This work was supported by the European Regional Development Fund (ERDF/ EFRE: ZS/2016/04/78156)
within the research center dynamic systems (CDS).
\bibliographystyle{siam}
\bibliography{references}

\end{document}